\documentclass[10pt, conference]{IEEEtran}
\ifCLASSINFOpdf
\else
\fi
\usepackage{amsmath,amssymb,amsfonts}
\usepackage{graphicx}
\usepackage{tabularx}
\usepackage{threeparttable}
\usepackage{multirow}
\usepackage{color}
\usepackage{ulem}

\newcommand{\ignore}[1]{}

\newif\ifsubmit
\submittrue
\ifsubmit

    \newcommand{\mert}[1]{}
    \newcommand{\carl}[1]{}
    \newcommand{\izzat}[1]{}
    \newcommand{\wengcho}[1]{}
    \newcommand{\levent}[1]{}
    \newcommand{\wenmei}[1]{}

    \newcommand{\todo}[1]{}
    \newcommand{\tocite}[1]{}
    \newcommand{\outline}[1]{}
    \newcommand{\reviewer}[2]{}
    
    \newcommand{\rmtext}[1]{}

\else

    \definecolor{commentColor}{rgb}{0.1, 0.6, 1.0}
    \definecolor{outlineColor}{rgb}{0.0, 0.50, 0.0}
    \definecolor{toCiteColor}{rgb}{1.0, 0.0, 0.5}
    \definecolor{addedColor}{rgb}{0.3, 0.5, 1.0}
    \definecolor{removedColor}{rgb}{0.8, 0.8, 0.8}

    \newcommand{\mert}[1]{[\footnote{{\color{commentColor}MH: #1}}]}
    \newcommand{\carl}[1]{[\footnote{{\color{commentColor}CP: #1}}]}
    \newcommand{\izzat}[1]{[\footnote{{\color{commentColor}IE: #1}}]}
    \newcommand{\wengcho}[1]{[\footnote{{\color{commentColor}WC: #1}}]}
    \newcommand{\levent}[1]{[\footnote{{\color{commentColor}LG: #1}}]}
    \newcommand{\wenmei}[1]{[\footnote{{\color{commentColor}WMH: #1}}]}

    \newcommand{\todo}[1]{[{\color{red}TODO: #1}]}
    \newcommand{\tocite}[1]{[{\color{toCiteColor}CITE: #1}]}
    \newcommand{\outline}[1]{[{\color{outlineColor}OUTLINE: #1}]}
    \newcommand{\reviewer}[2]{[\footnote{{\color{red}Reviewer: #1} {\color{outlineColor}Response: #2}}]}
    
    \newcommand{\rmtext}[1]{{\color{removedColor}\sout{#1}}}

\fi

\begin{document}
%
\title{Fast Numerical Integration Techniques for 2.5-Dimensional Inverse Problems}




%
\author{\IEEEauthorblockN{Mert Hidayeto\u{g}lu, Michael Oelze, Erhan Kudeki, and Weng Cho Chew}
\IEEEauthorblockA{Department of Electrical and Computer Engineering\\
University of Illinois at Urbana-Champaign,
Urbana, IL 61801, USA}}


\maketitle


%

\begin{abstract}
Inverse scattering involving microwave and ultrasound waves require numerical solution of nonlinear optimization problem. To alleviate the computational burden of a full three-dimensional (3-D) inverse problem, it is a common practice to approximate the object as two-dimensional (2-D) and treat the transmitter and receiver sensors as 3-D, through a Fourier integration of 2-D modes of scattering. The resulting integral is singular, and hence requires a prohibitively large number of integration points, where each point corresponds to a 2-D solution. To reduce the computational complexity, this paper proposes fast integration approaches by a set of transformations. We model the object in 2-D but the transmit and receiver pairs as 3-D; hence, we term the solution as a 2.5-D inverse problem. Convergence results indicate that the proposed integration techniques have exponential convergence and hence have a reduces the computational complexity to compute 2.5-D Green's function to solve inverse scattering problems.
\end{abstract}

\section{Introduction}

Wave scattering involving inhomogeneous media requires numerical solutions of the wave equation in a three-dimensional (3-D) volume. Such solutions have a plethora of applications in inverse problems such as ultrasound and microwave imaging \cite{lavarello2009tomographic,neira2017high,catapano20093d, zakaria2010finite,wiskin2020full, hesford2010fast}. 
Volume integral equation (VIE) formulation provides stable numerical solutions via time-harmonic fields. In this case, the volume is discretized using basis and testing functions, and the unknown coefficients of the basis functions are found by using the method of moments~\cite{harrington1993field} through matrix inversion. Nevertheless, solving VIE with a higher dimension can be computationally costly~\cite{bellman1966dynamic}. That is, solving the 3-D problem involves larger number of unknowns by an order of magnitude compared to solving a 2-D problem (which can also be costly at large scale~\cite{hidayetoglu2018fast}). In many scenarios~\cite{nikolova2011microwave}, the object to be imaged can be well-approximated as 2-D; however, the finite dimensions of physical transmitters (Tx) and receivers (Rx) require modelling them in 3-D.

\begin{figure}[t]
  \centering\includegraphics[width=0.48\textwidth]{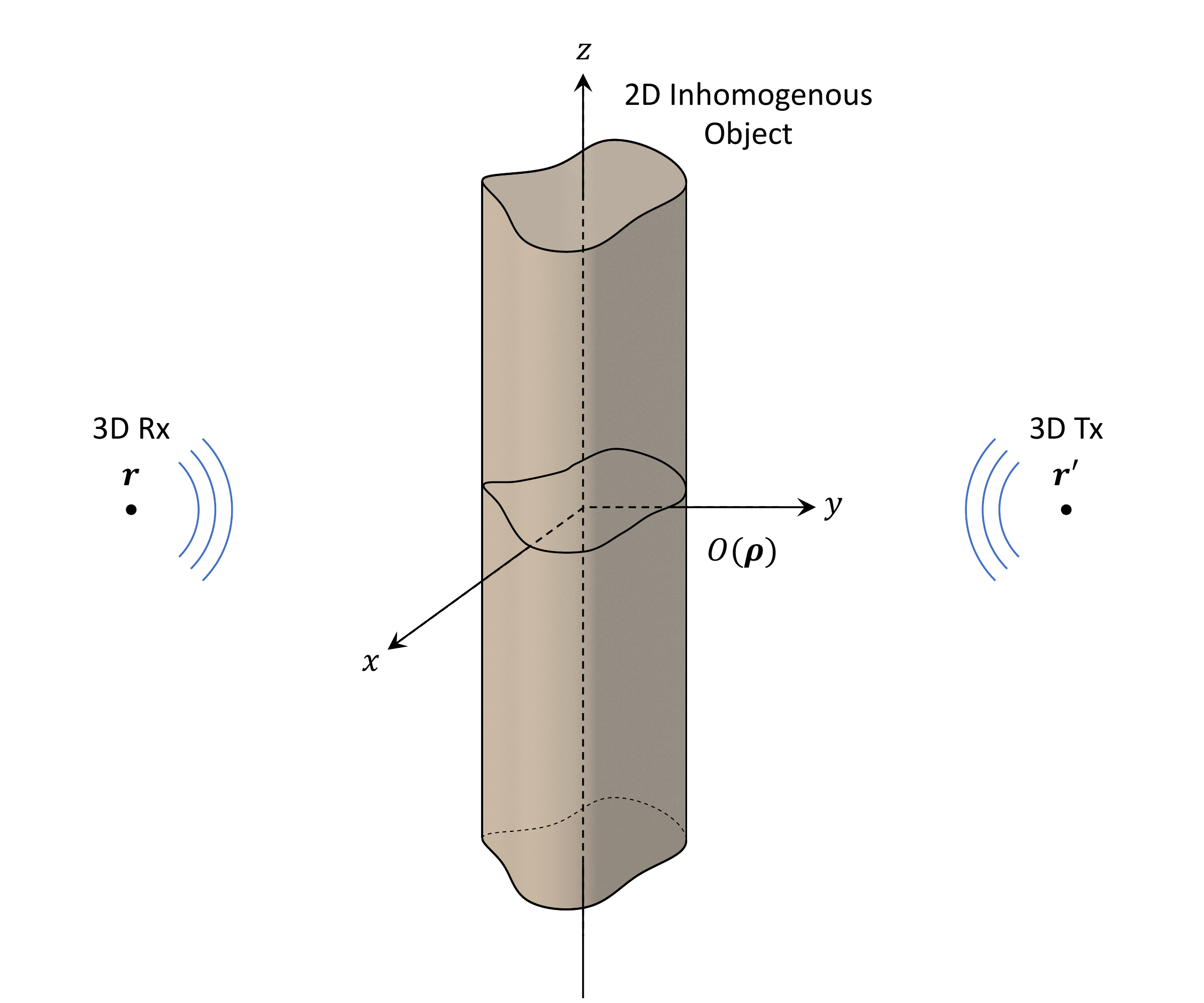}
  \caption{In 2.5-D inverse problem, the inhomogeneous object is modeled in 2-D, and the Tx/Tr pairs are modeled in 3-D.}\label{fig:2p5D_problem}
\end{figure}

To reduce the dimensionality, 2.5-D formulation models the unknown inhomogeneous scatterer as 2-D, but the Tx/Rx are modeled as 3-D, (see Fig.~\ref{fig:2p5D_problem}). To that extent, the 2-D modes of radiation from a 3-D Tx can be computed through the 2.5-D Green's function. In that case, the incident field is expanded with cylindrical Fourier modes so that the scattering with each mode can be solved in an embarrassingly parallel manner, i.e., each processor solves one mode of scattering. Next, the scattered field must be integrated back with inverse Fourier synthesis to find the measured field at the 3-D Rx.  Nevertheless, the synthesis with the 2.5-D Green's function is not numerically attainable because of two major drawbacks related to numerical implementation with quadrature rules: 1) the integral has infinite limits and 2) the integral involves a logarithmic singularity, prohibiting the convergence of the integral.

In this paper, we analyze the 2.5-D Green's function and its singularity. For improving convergence of numerical integrals, we follow a set of analytical transformations~\cite{brekhovskikh2013acoustics} in the complex plane to suppress the singularity and use the steepest-descent path to provide for exponential convergence~\cite{chew1990waves}. Moreover, we propose parametrizations for mapping infinite integral path to a finite interval for numerical integration. Numerical results demonstrate exponential convergence when integrating the 2.5-D Green's function to obtain 3-D radiation from a point source. As a minor contribution, we propose a set of regularization techniques to handle pathological cases. The results of this paper will yield efficient numerical formulation to solve 2.5-D inverse scattering and imaging problems.

The rest of the paper is organized as follows: Section II introduces the 2.5-D Green's function for scalar fields. Section III analyzes the branch cuts and the singularity in the complex plane that restricts the convergence of numerical integration. Section IV provides an angular spectrum transformation to suppress the singularity of the 2.5-D Green's function. Section V presents the proposed techniques for application on 2.5-D inverse problems. Section VI presents numerical convergence results of the proposed integration techniques. Finally, Section VII discusses the results and conclusions of the study.

\section{2.5-D Green's Function}

With $e^{-i\omega t}$ time convention, the three-dimensional (3-D) Green's function of the wave equation for a homogeneous medium can be defined as
\begin{equation}\label{eq:3D_green}
    g_0(\boldsymbol{r}, \boldsymbol{r}') = \frac{e^{ik_0|\boldsymbol{r}-\boldsymbol{r}'|}}{4\pi|\boldsymbol{r}-\boldsymbol{r}'|},
\end{equation}
which is the solution of the wave equation
\begin{equation}\label{eq:wave_freq}
    \nabla^2g_0(\boldsymbol{r},\boldsymbol{r}')+k_0^2g_0(\boldsymbol{r},\boldsymbol{r}')=-\delta(\boldsymbol{r}-\boldsymbol{r}'),
\end{equation}
with a complex wavenumber, representing a lossy medium. Here, the right-hand-side of (\ref{eq:wave_freq}) describes a point source located at $\boldsymbol{r}'=\hat{\boldsymbol{x}}x'+\hat{\boldsymbol{y}}y'+\hat{\boldsymbol{z}}z'$ in the Cartesian coordinates and $\delta(\boldsymbol{r})$ is a shorthand for $\delta(x)\delta(y)\delta(y)$, where $\delta(x)$ is the Dirac-delta function.

The 3-D Green's function in (\ref{eq:3D_green}) can be expanded by cylindrical harmonics \cite{chew1990waves} as
\begin{equation}\label{eq:fourier_synthesis}
	g_0(\boldsymbol{r},\boldsymbol{r}') = \frac{i}{8\pi}\int_{-\infty}^\infty dk_z\,H_0^{(1)}(k_\rho|\boldsymbol{\rho}-\boldsymbol{\rho}'|)e^{ik_z(z-z')}
\end{equation}
where $\boldsymbol{\rho}=\boldsymbol{\hat{x}}x+\boldsymbol{\hat{y}}y$ represents the two-dimensional (2-D) coordinate vector and $H_0^{(1)}$ is the zeroth-order Hankel function of the first kind. Here, $k_\rho = \sqrt{k_0^2-k_z^2}$ and can take complex values, depending on the values of $k_0$ and $k_z$. As a result of the Fourier expansion in (\ref{eq:fourier_synthesis}), we can deduce that
\begin{equation}\label{eq:2p5_green}
	g_0(\boldsymbol{\rho},\boldsymbol{r}') = \frac{i}{4}H_0^{(1)}(k_\rho|\boldsymbol{\rho}-\boldsymbol{\rho}'|)e^{-ik_zz'}
\end{equation}
is the spatial Fourier transform of the 3-D Green's function in the $z$ direction and describes a $k_z$ mode of radiation from a point source at $\boldsymbol{r}'$ in 3-D. In this paper, we define~(\ref{eq:2p5_green}) as the 2.5-D free-space Green’s function and the observation point is denoted with $\boldsymbol{\rho}$ in 2-D, and the $k_z$ is the Fourier transform variable and for brevity, the Fourier mode is implied. For example, the conventional 2-D Green's function can be defined as $k_z=0$ mode of (\ref{eq:2p5_green}), with no $z$ variation in the 3-D space. The derivation and interpretation of the 2.5-D Greens function is provided in the appendix (\S\ref{sec:appendix}).



\begin{figure}[t]
  \centering\includegraphics[width=0.45\textwidth]{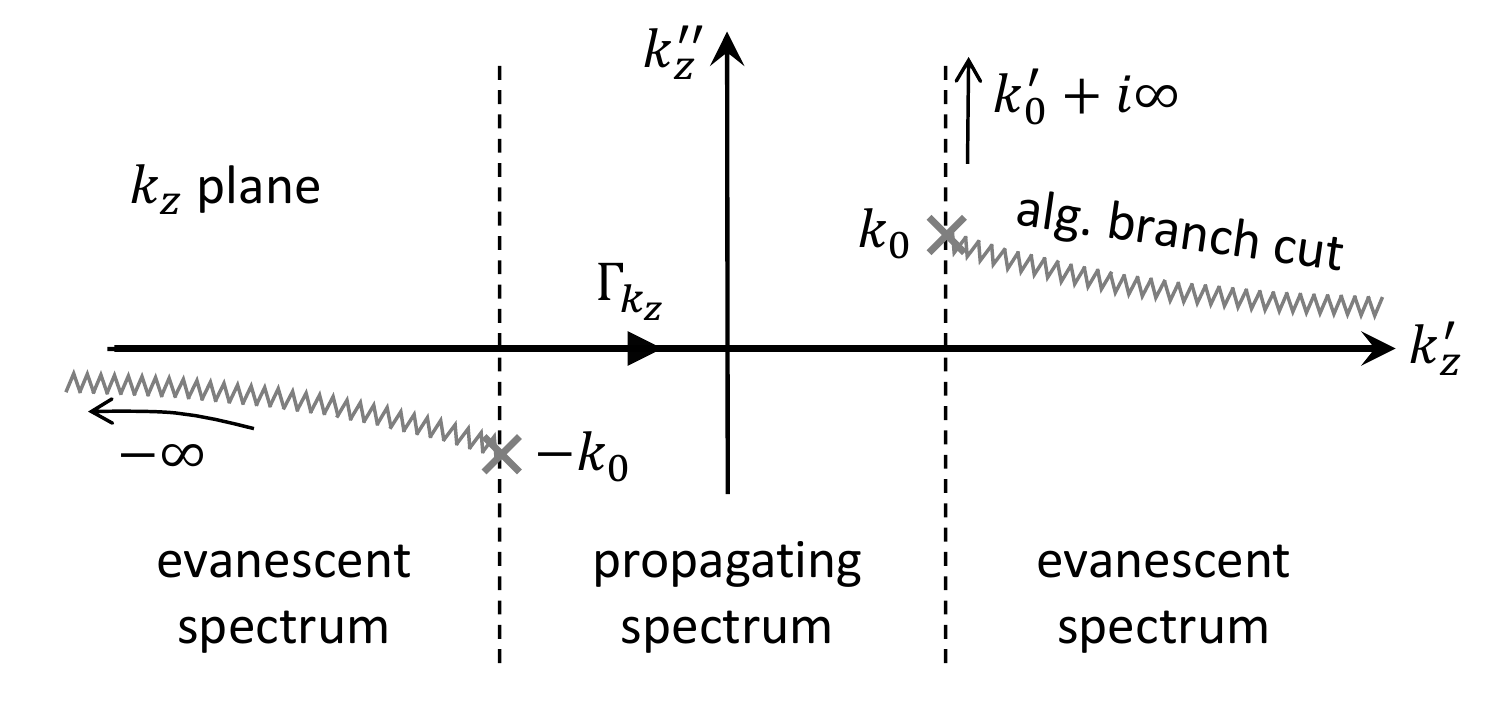}
  \caption{Algebraic branch cuts in the complex $k_z$ plane. The branch points at $\pm k_0$ shift to the real line when the medium is lossless. The evanescent and propagating $k_z$ spectrum are approximately separated with the dashed lines.}\label{fig:kz_plane}
\end{figure}

\begin{figure}[t]
  \centering\includegraphics[width=0.45\textwidth]{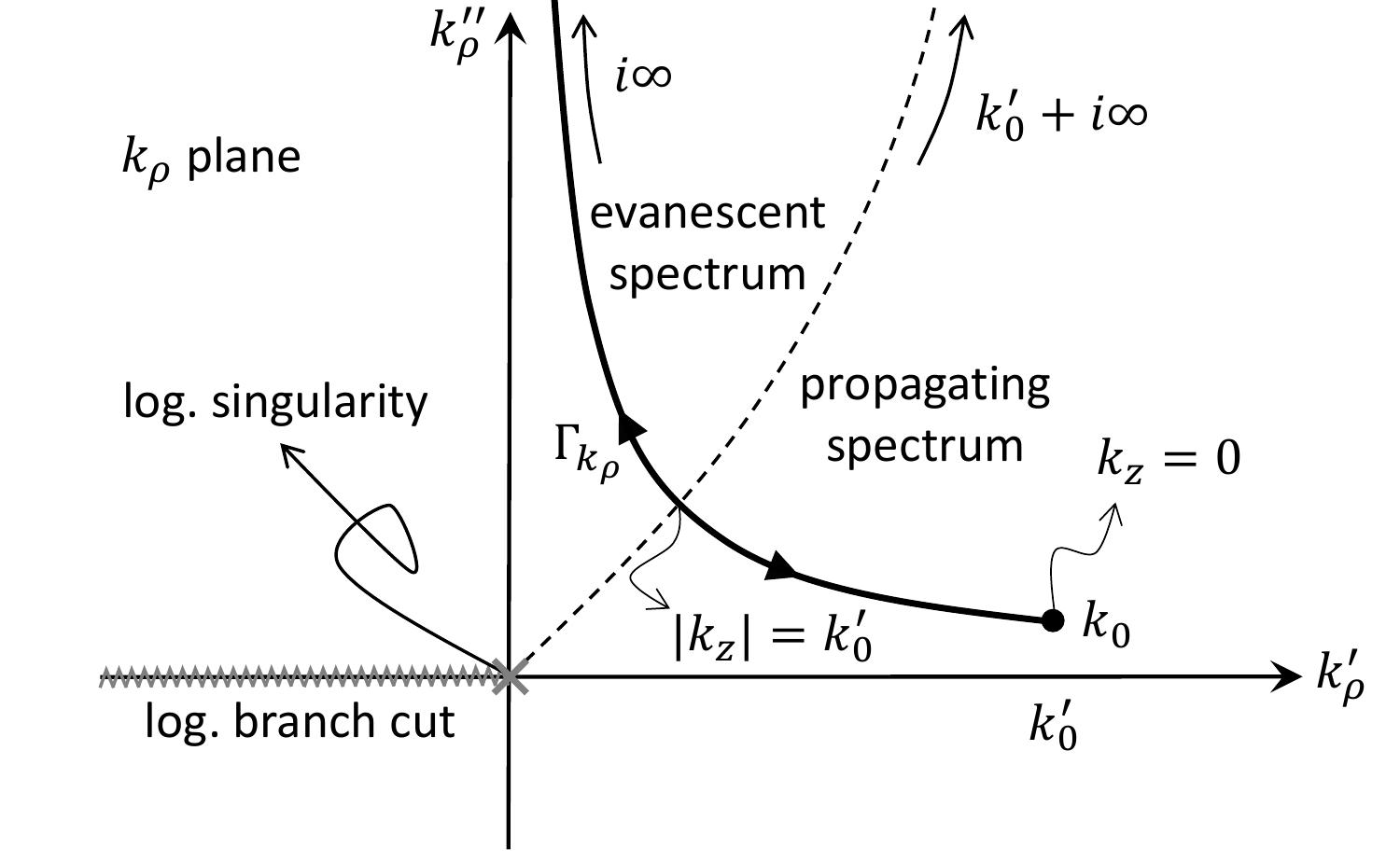}
  \caption{Logarithmic branch cut in the $k_z$ complex plane. When the medium is lossless, the integration path $\Gamma_{k_\rho}$ approaches to the $x$ and $y$ axes, and to the logaritmic singularity at the origin.}\label{fig:kr_plane}
\end{figure}

\section{Singularity of 2.5-D Green's Function}\label{sec:formulation_rectangular}

The Fourier synthesis in (\ref{eq:fourier_synthesis}) involves a logarithmic singularities at $k_z=\pm k_0$, when $k_\rho=0$. That spoils the numerical integration of the 2.5-D Green's function. For analyzing the algebraic branch points and logarithmic singularities, we simplify (\ref{eq:fourier_synthesis}) by placing the source at the origin, i.e.,
\begin{equation}\label{eq:hankel}
    g_0(r)=\frac{i}{8\pi}\int_{-\infty}^\infty dk_z \,H_0^{(1)}(k_\rho\,\rho)e^{ik_zz},
\end{equation}
without any loss of generality. In the integrand of (\ref{eq:hankel}), the exponential function in the second term is entire (analytic everywhere), and therefore has no complication. However, the first term, the Hankel function is a nested function. The first function maps $k_z$ to $k_\rho = \sqrt{k_0{}^2-k_z^2}$. The second function maps $k_\rho$ to $H_0{}^{(1)}(k_\rho\,\rho)$. Sections \S~\ref{sec:kz_plane} and \S~\ref{sec:kr_plane} describes critical points on the $k_z$ and $k_\rho$ complex planes, respectively.

\subsection{Algebraic Branch Cuts in the $k_z$ Plane}\label{sec:kz_plane}
Figure~\ref{fig:kz_plane} shows the complex $k_z$ plane for a lossy medium, where $\Gamma_{k_z}$ is the integration contour on the real line. The $k_z$ plane involves two algebraic branch points
at $k_z=\pm k_0$ that exist with branch
cuts. The standard branch cuts are defined from branch points extending to infinity, parametrized as $k_z=\pm\sqrt{k_0^2+c^2}$ with a positive $c$. Here $k_0=k_0'+ik_0''$ is complex, where $k_0'$ and $k_0''$ are real quantities representing the spatial frequency and attenuation terms, respectively.

When the medium is lossless, the branch points shifts to the real axis, which coincides with $\Gamma_{k_z}$. In that case, $|k_z|\leq k_0$ corresponds to the propagating spectrum and the $|k_z|>k_0$ corresponds to the evanescent spectrum of the 3-D Green's function. When the medium is lossy, the branch points $k_0$ and $-k_0$ shifts to the upper- and lower-half of the $k_z$-plane, respectively. In this case, the propagating and evanescent spectra are not physically meaningful, nevertheless, the branch points follow the straight boundaries shown with the dashed line in Fig.~\ref{fig:kz_plane}.

\subsection{Logarithmic Singularity in the $k_\rho$ Plane}\label{sec:kr_plane}
Figure~\ref{fig:kr_plane} shows the complex $k_\rho$ plane for a lossy medium. Mapping from $\Gamma_{k_z}$ to $\Gamma_{k_\rho}$ on the complex $k_\rho$ plane is represented with bold line with arrows. The integration path $\Gamma_{k_\rho}$ draws the parabola $k_\rho = \sqrt{k_0^2-k_z^2}$. It is worthy to note that the root of the first quadrant of the complex $k_\rho$ plane is designated as the physical root. When the medium is lossless, that parabola becomes a rectangular path coincides with the $x$ and $y$ axes, as well as the singularity at the origin. Fig.~\ref{fig:grid} in the appendix shows the loci of $k_\rho$ for various complex values of $k_0$.

The complex plane plane in Fig.~\ref{fig:kr_plane} involves logarithmic singularities at the origin and at the infinity. These singularities are the endpoints of the branch cut, defined on the negative real line. The dashed line depicts the boundary between respective evanescent and propagating $k_z$ spectra in the lossless case. As evident from the analysis in this section, the $k_z$ integration contour intersects with the logarithmic singularity at $k_z=\pm k_0$ in lossless media. That hampers the convergence of the Fourier synthesis in (\ref{eq:hankel}).

\section{Angular Spectrum Formulation}\label{sec:formuation_angular}
A faster convergence of (\ref{eq:hankel}) can be obtained by the angular-spectrum formulation~\cite{brekhovskikh2012waves}. This removes the algebraic branch points from the integral by a mapping from the $k_z$ spectrum into the angular spectrum, i.e., by a change of variable $k_z=k_0\sin\theta$ and $k_\rho=k_0\cos\theta$. Then
\begin{equation}
    g_0(r) = \frac{i}{8\pi}\int_{\Gamma_{\theta}}d\theta\,k_0\cos\theta\, H_0^{(1)}(k_0\cos\theta\,\rho)e^{ik_0\sin\theta\,z}.\label{eq:angular}
\end{equation}
In this mapping, $\theta$ takes complex values when $|k_z|>k_0'$, i.e., $\theta=\theta'+i\theta''$. We can parameterize $\Gamma_\theta$ in three partitions, where the real and complex partitions involve the propagating and evanescent $k_z$ spectra, respectively.

\subsection{Logarithmic Singularities in the $\theta$ Plane}

\begin{figure}[t]
    \centering\includegraphics[width=0.5\textwidth]{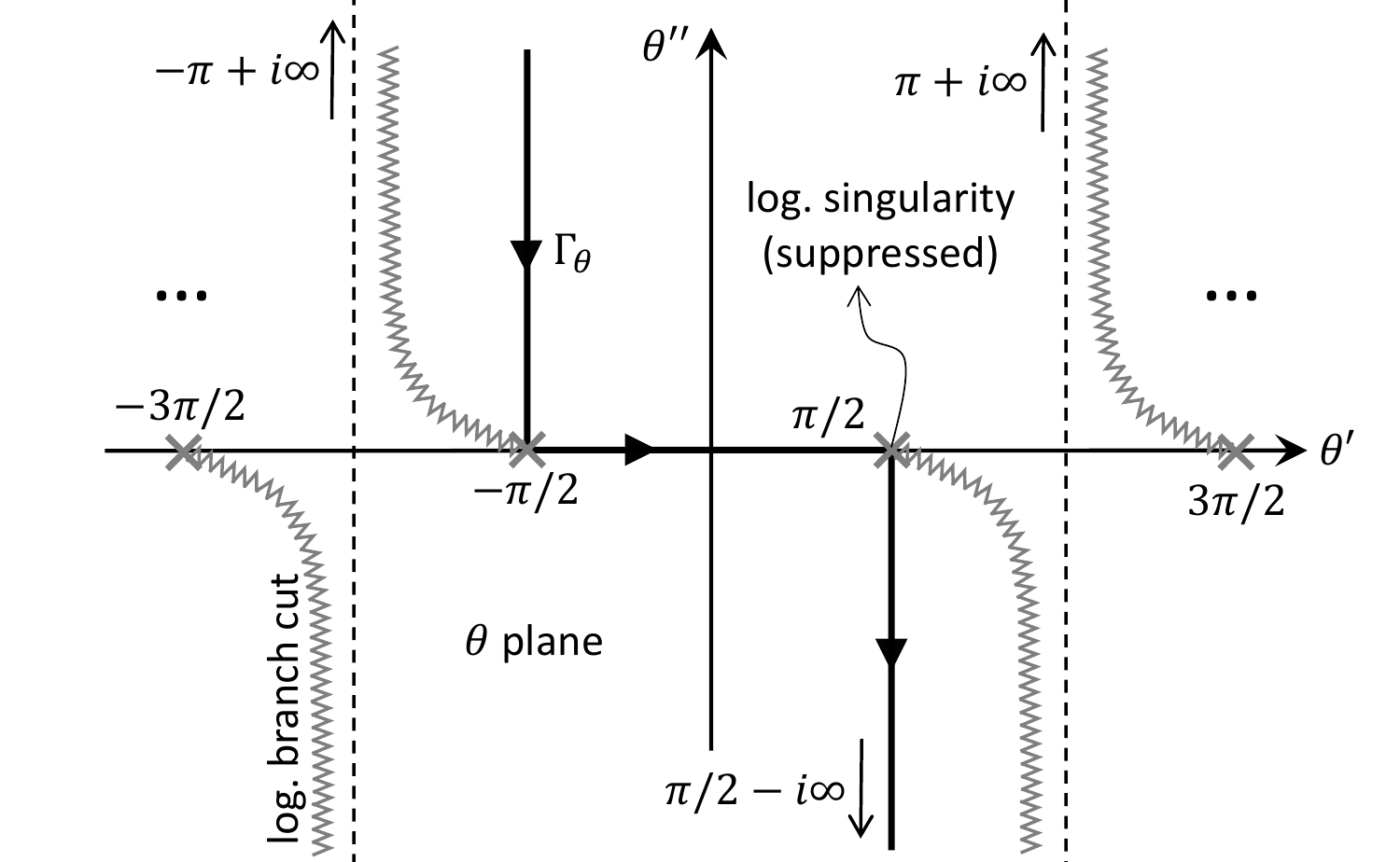}
  \vspace{-3mm}
  \caption{Angular integration path $\Gamma_\theta$ on the $\theta$ plane which has suppressed logarithmic singularities and removes algebraic branch points. Real and imaginary $\theta$ involves propagating and evanesecent spectra, respectively.}\label{fig:kt_plane}
\end{figure}

Fig.~\ref{fig:kt_plane} shows the integration path $\Gamma_{\theta}$ on the $\theta$ plane. The figure also depicts the the 2$\pi$-periodic logarithmic branch cuts parametrized standardly as $\theta=\cos^{-1}(-c/k_0)$ and logarithmic singularities at $\theta=\pm\pi/2+2\pi p$, where $p$ is an integer, as a result of the Hankel function with a trigonometric argument in (\ref{eq:angular}). The branch cuts take rectangular shapes when the medium is lossless. Note that the logarithmic singularities are suppressed by the first term of the integrand. We can observe that by recalling the asymptotic form of the Hankel function for small argument, i.e., $H_0^{(1)}\sim \log$, and letting
\begin{equation*}
    \lim_{\theta\rightarrow\pi/2}\cos\theta\log(k_0\cos\theta\,\rho)=0
\end{equation*}
The very same fact prevails at $\theta=-\pi/2$ as well as at other logarithmic singularities in the $\theta$ plane. 

Supressing the logarithmic singularities with the angular formulation does not change the fact the integrand is not analytic at the logarithmic branch points $\theta=\pi/2+p\pi$, which prevents numerical integrations. Notwithstanding, the suppression reduces the integration error and helps convergence as evident in numerical results because the integrand, and therefore the integration error, is small around the branch points. Note that imposing a loss does not help convergence in this formulation, as opposed to the original formulation in (\ref{eq:hankel}) of the previous section because the branch points and integration contour remains fixed regardless of the loss term as evident in numerical results (Fig.~\ref{fig:result_1}).

\section{Fast Integration Approaches}
By analytical continuity, we can deform the spectral integration path so as to avoid rapidly varying and unstable domains of the oscillatory integrand in (\ref{eq:angular}). In this spirit, we propose two integration approaches; the first one follows the stationary-phase point, and the second one follows the steepest-descent path.

The spectral integration (\ref{eq:angular}) can be written as a multiplication of slowly- and rapidly-varying parts as
\begin{equation}\label{eq:parts}
    g(r) = \frac{i}{8\pi} \int_{\Gamma_\theta} d\theta\,\widehat{H}(\theta) e^{f(\theta)},
\end{equation}
where 
\begin{equation}
    \widehat{H}(\theta) = \frac{k_0\cos\theta H_0^{(1)}(k_0\cos\theta\rho)}{e^{ik_0\cos\theta\rho}}
\end{equation}
is the slowly-varying part and $e^{f(\theta)}$ is the rapidly-varying part. The integrand is highly-susceptible to change in $\theta$ in the exponent of the rapidly-warying part
\begin{equation}\label{eq:exponent}
    f(\theta) = ik_0r\cos(\theta-\theta_0),
\end{equation}
where $\theta_0=\tan^{-1}(z/\rho)$ is the observation angle (which is always real) and $r=\sqrt{\rho^2+z^2}$ is the observation distance from the source. It is also worthy to note that the slowly-varying part is singular when $\theta_0 = \pm\pi/2$ since $\rho=0$.

\begin{figure}[t]
    \centering\includegraphics[width=0.5\textwidth]{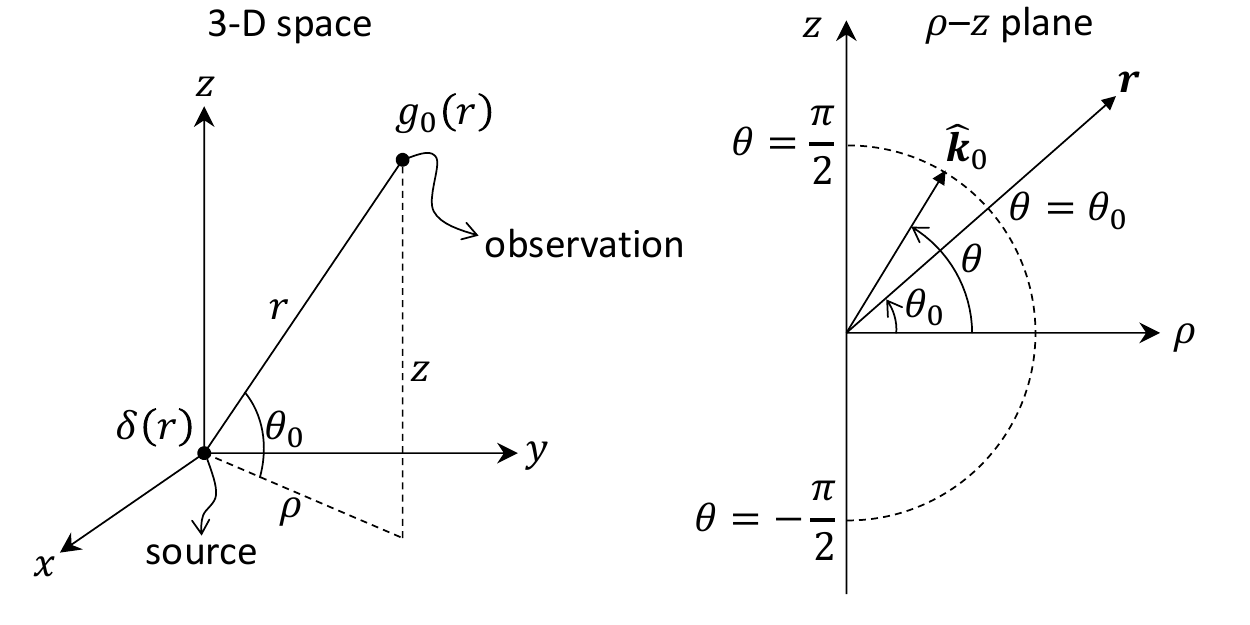}
  \vspace{-3mm}
  \caption{A point source at the origin and observation at distance $r$. The angle of the observation point to the $x$-$y$ plane corresponds to the saddle point in the angular spectrum formulation.}\label{fig:observation}
\end{figure}

\begin{figure}[t]
    \hspace{-3mm}\centering\includegraphics[width=0.49\textwidth]{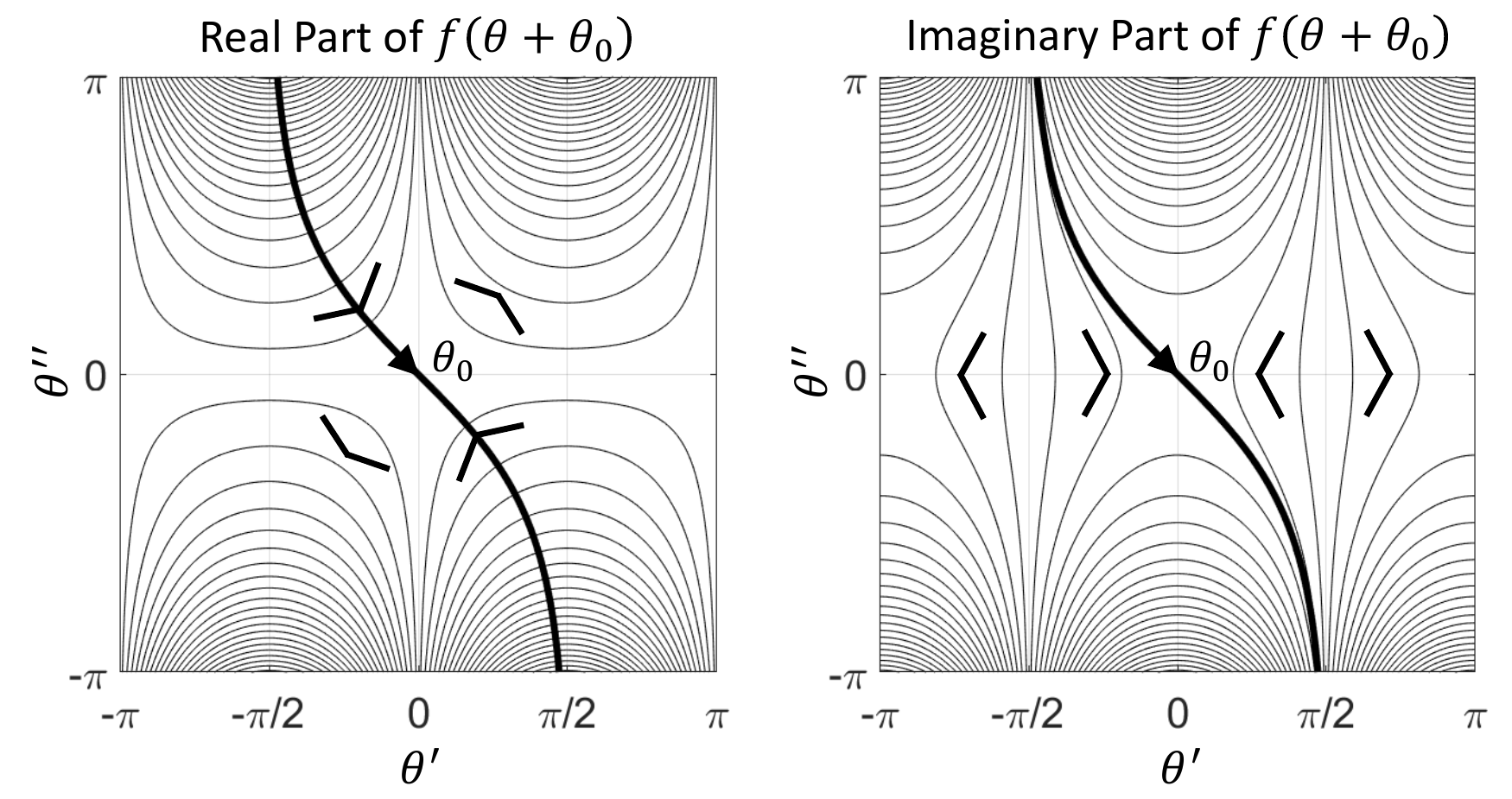}\\
  {\hspace{0.25cm}\small(a)\hspace{4.25cm}(b)}\\
  \vspace{4mm}
  \hspace{-3mm}\centering\includegraphics[width=0.49\textwidth]{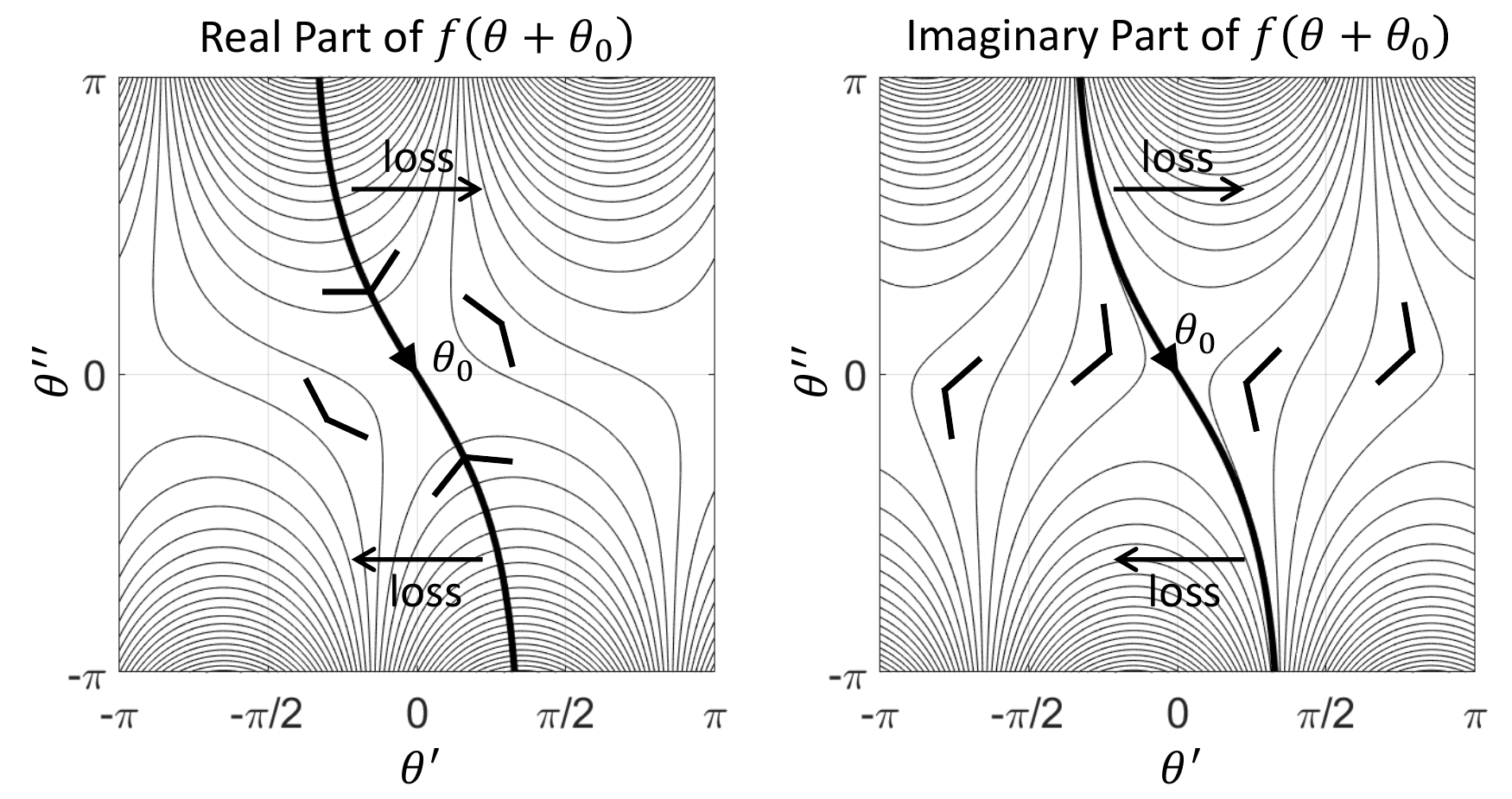}\\
  {\hspace{0.25cm}\small(c)\hspace{4.25cm}(d)}
  \caption{Contour plots of (a),(c) real  and (b),(d) imaginary parts of $f(\theta+\theta_0)$ and the steepest-descent paths for (a),(b) lossless ($k_0=2\pi$) and (b),(d) lossy ($k_0=2\pi[1+i0.5]$) medium. The ascent directions are depicted with large arrow heads.}\label{fig:integral_paths}
\end{figure}

By the maximum modulus principle~\cite{palka1991introduction}, the analytic function $e^{f(\theta)}$ cannot contain a local maximum or minimum anywhere in the complex plane, instead, it has saddle points, which  physically correspond to the stationary phase points. Fig.~\ref{fig:integral_paths} shows the controur plots of real and imaginary parts of $f(\theta)$, where the saddle point is at the origin. The idea of the stationary-phase approximation is to choose an integral path so as to pass through the saddle point, align it in the direction of fastest decay, and truncate it as it captures the important information. Investigating the derivative \mbox{$f'(\theta) = -ik_0r\sin(\theta-\theta_0)$},
we can deduce that the saddle point is at $\theta_0$ since $f'(\theta_0)=0$. This saddle point corresponds to the observation angle as depicted in Fig.~\ref{fig:observation}.

\subsection{Stationary-Point Integration: Quadratic Approximation}\label{sec:quadratic_approx}
The stationary-phase point contributes to the integral the most, and therefore we want to modify the integration path so that it passes through $\theta_0$ by contour deformation. Consequently, we can simply truncate the integral around $\theta_0$. To do that, we can expand $f(\theta)$ around $\theta_0$ as
\begin{equation}
    f(\theta) = f(\theta_0) - \frac{f^{(2)}(\theta_0)}{2!}(\theta-\theta_0)^2  + \frac{f^{(4)}(\theta_0)}{4!}(\theta-\theta_0)^4 + \dots
\end{equation}
and letting $\theta=\theta_0+te^{i\hat{\theta}}, t\ge 0$, we can approximate
\begin{equation}
    f(\theta)-f(\theta_0)\approx-\frac{ik_0rt^2}{2}e^{i2\hat{\theta}}
\end{equation}
with an error term of $\mathcal{O}(t^4)$. The steepest descent path lies in the constant imaginary part of $f(\theta)$, therefore by setting $\Im\{f(\theta)-f(\theta_0)\}=0$, i.e., $\cos(2\hat{\theta}+\alpha)=0$, the steepest directions at $\theta_0$ can be found as
\begin{equation}
    \hat{\theta} = \frac{\pi/2+\pi p-\alpha}{2},
\end{equation}
where $p$ is an integer and $\alpha$ is the loss angle, i.e., $k_0=|k_0|e^{i\alpha}$. As a result, there are four steepest directions. Investigating the real part in the steepest directions,
\begin{equation}
\Re\{f(\theta)-f(\theta_0)\}=\frac{rt^2}{2}\sin(\pi/2+\pi p)
\end{equation}
we can see that with increasing $t$, the exponent increases for an odd $p$ and decreases for an even $p$. Observing the downward direction of original path shown in Fig.~\ref{fig:kt_plane}, we can choose $p=-3$ and write the linear fit to the steepest-descent path as
\begin{equation}
    \theta = \theta_0 + te^{-i[\pi/4+\alpha/2]}\label{eq:parameter_approximate}
\end{equation}
for increasing $t$. By this parametrization, the angular integral can be approximated as
\begin{equation}
   g_0(r)\approx\frac{i}{8\pi}\int_{-\epsilon}^{\epsilon} dt\, \hat{H}(\theta)e^{i\boldsymbol{k}_0(\theta)\cdot\boldsymbol{r}}e^{-i[\pi/4+\alpha/2]},\label{eq:approximate_descent}
\end{equation}
where $\epsilon$ is a truncation parameter. When $k_0r$ is large, even a small $\epsilon$ yields sufficiently accurate integral results. In contrast, when $k_0r$ is small, we need to increase $\epsilon$ for good accuracy. However, because the quadratic approximation is accurate only around $\theta_0$, $\epsilon$ should not be selected too large or the integration path may cross a branch cut or the integrand goes into unstable areas in the angular spectrum. The latter case is shown in numerical results (Fig.~\ref{fig:result_2}).


In order to identify unstable areas in the $\theta$ plane, we can write $\theta=\theta'+i\theta''$ and $k_0=|k_0|e^{i\alpha}$ to expand (\ref{eq:exponent}) as
\begin{align*}
    f(\theta_0+\theta) &= ik_0r\cos\theta\\
    &=\frac{r|k_0|}{2}\left[e^{-\theta''+i(\pi/2+\alpha+\theta')}+e^{\theta''+i(\pi/2+\alpha-\theta')}\right].
\end{align*}
We can find the magnitude by considering the real part
\begin{align*}
    \Re\{f(\theta_0+\theta)\}= r|k_0|e^{-\theta''}\sin&(\alpha+\theta')/2\\
   & +r|k_0|e^{\theta''}\sin(\alpha-\theta')/2.
\end{align*}
Because the real part of the exponent grows and shrinks exponentially when traveling away from the real line, the exponential term also grows or shrinks exponentially, depending on the loss ($\alpha$) and the real angle ($\theta')$. The growing and shrinking areas for the lossless medium can be seen in Fig.~\ref{fig:kt_plane_efficient}. The stationary-phase point is located at a saddle point. An efficient integration path should avoid the growing areas, shaded with gray in the figure.

We can reparametrize (\ref{eq:parameter_approximate}) by writing $\theta=\theta'+i\theta''$ and
\begin{equation}
    \theta'' = \theta_0+\theta'[1-i\tan(-\pi/4-\alpha/2)],
\end{equation}
and for lossless medium, $\alpha=0$ and we can safely set $\epsilon=\pi/2$ (one can geometrically confirm this from Fig.~\ref{fig:kt_plane_efficient}), yielding
\begin{equation}
    g_0(\boldsymbol{r},\boldsymbol{r}')\approx\frac{1}{2\pi}\int_{-\pi/2}^{\pi/2} d\theta'\,[1-i]k_\rho g(\boldsymbol{\rho},\boldsymbol{r}')e^{ik_zz},
\end{equation}
where
\begin{equation}
    g_0(\boldsymbol{\rho},\boldsymbol{r}')=\frac{i}{4}H_0^{(1)}(k_\rho|\boldsymbol{\rho}-\boldsymbol{\rho}'|)e^{-ik_zz'}
\end{equation}
is the 2.5-D Green's function, $k_\rho = k_0\cos\theta$, and $k_z = k_0\sin\theta$.

\subsection{Steepest Descent Integration on the $\theta$ Plane}\label{sec:exact_SD}
Since $f(\theta)$ is analytic, it satisfies the Cauchy-Riemann equations~\cite{palka1991introduction}, i.e.,
\begin{equation}
    \nabla\Re\{f(\theta)\}\cdot\nabla\Im\{f(\theta)\} = 0.
\end{equation}
This implies that the contour of the steepest descent is also the contour of constant phase, i.e., $\nabla\Im\{f(\theta)\}=0$. We follow the constant phase by the saddle point by equating \begin{equation}\label{eq:constant_phase}
    \Im\{f(\theta_0+\theta)\}=\Im\{f(\theta_0)\}.
\end{equation}
Then we can parametrize the steepest-descent contour as
\begin{equation}\label{eq:imaginary}
    \cos\theta'\cosh\theta'' + \tan(\alpha)\sin\theta'\sinh\theta''= 1
\end{equation}
and therefore $\theta = \theta_0+\theta'+i\theta''$. We can solve the above in a closed form as
\begin{equation}
    \theta''= \log\left(\frac{\cos\alpha-\sin\theta'}{\cos(\alpha-\theta')}\right)
\end{equation}
as a function of
\begin{equation*}
    -\pi/2+\alpha < \theta' < \pi/2-\alpha.
\end{equation*}

\begin{figure}[t]
    \centering\includegraphics[width=0.45\textwidth]{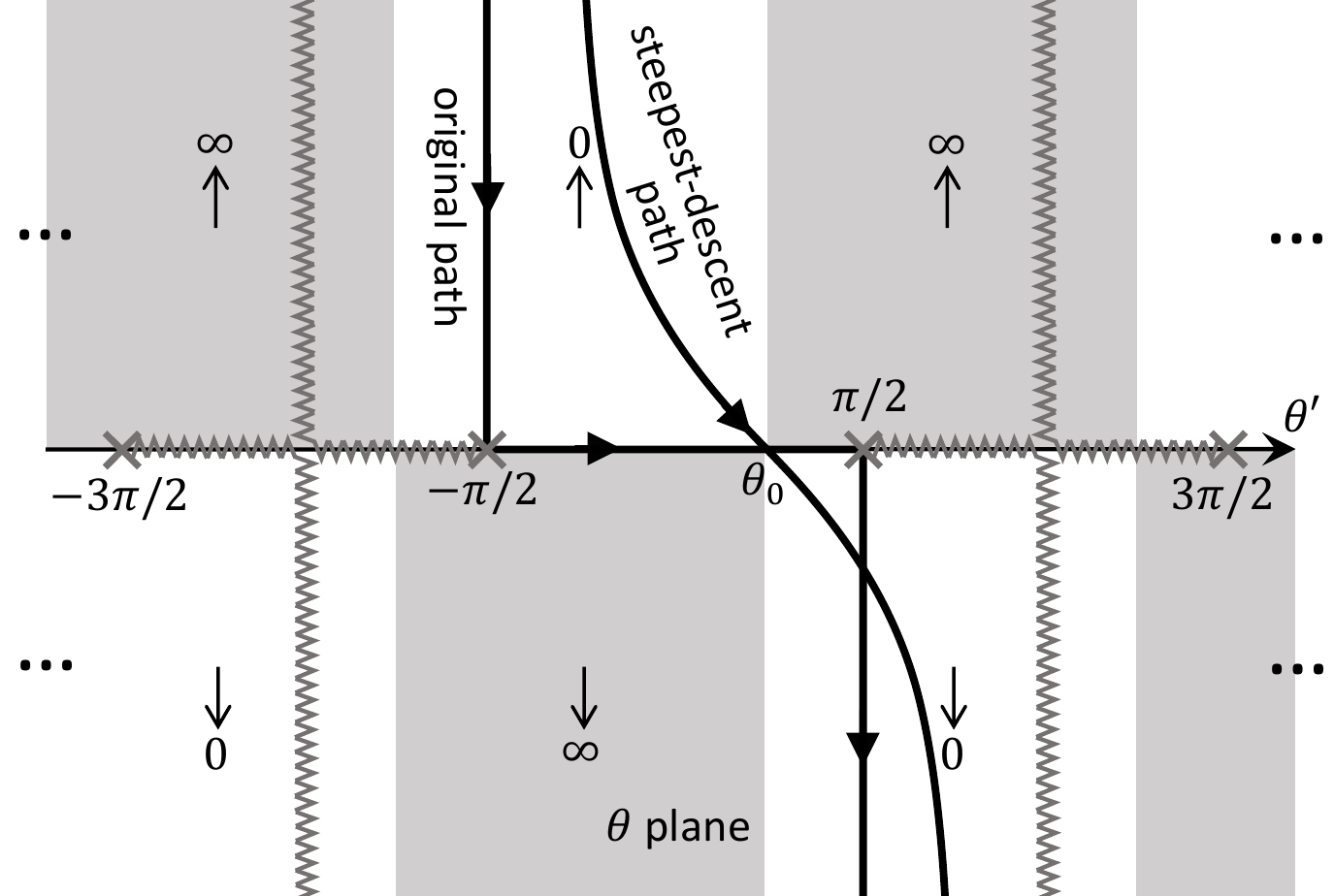}
  \caption{Angular and steepest-descent integration paths in the complex $\theta$ plane for lossless medium. The darker shades depicts unstable regions where the integrand grows exponentially.}\label{fig:kt_plane_efficient}
\end{figure}

Fig.~\ref{fig:integral_paths} shows the steepest-descent path, where the imaginary part of the exponent stays constant. Also, Fig.~\ref{fig:kt_plane_efficient}, we can see that the steepest-descent and original integrals are equal by contour deformation. By noting $d\theta = c(\theta')d\theta'$, where
\begin{equation}
    c(\theta') = 1 - i\left[\frac{\cos(\theta')}{\cos\alpha-\sin\theta'}+\tan(\alpha-\theta')\right]
\end{equation}
is the correction term for the transformation, the steepest-descent integral can be explicitly written as
\begin{equation}
    g_0(r) = \frac{i}{8\pi}\hspace{-2mm}\int\limits_{-\pi/2+\alpha}^{\pi/2-\alpha}\hspace{-2mm} d\theta' c(\theta')k_0\cos\theta H_0^{(1)}(k_0\cos\theta\rho)e^{ik_0\sin\theta z}.\label{eq:steepest_descent}
\end{equation}

Note that, we converted the infinite integral in (\ref{eq:fourier_synthesis}) exactly to a finite steepest-descent integral in (\ref{eq:steepest_descent}). For numerical implementation of the above with Gaussian-Legendre (G-L) quadrature, we propose an additional transformation
\begin{equation}
    t=\frac{\theta'}{\pi/2-\alpha}
\end{equation}
to standardize the expression. As a result, we can use a legacy table of G-L quadrature points and weights to integrate
\begin{equation}
    g_0(\boldsymbol{r},\boldsymbol{r}') = \frac{1}{2\pi}\int_{-1}^1dt\,[\pi/2-\alpha] c(\theta')k_\rho g_0(\boldsymbol{\rho},\boldsymbol{r}') e^{ik_z z},\label{eq:3greens_function_eff}
\end{equation}
where $g_0(\boldsymbol{\rho},\boldsymbol{r}')$ is the 2.5-D Green's function given in (\ref{eq:2p5_green}). When implementing the G-L quadrature within the \texttt{for} loop, the following variables are computed in the order of their dependencies: $t\rightarrow\theta'\rightarrow\theta''\rightarrow\theta\rightarrow k_\rho\,,k_z$. As a drawback of the mapping of the infinite interval $k_z\in(-\infty,\infty)$ into a finite interval $t\in[-1,1]$, the integrand is singular at the limits of the interval and therefore cannot be evaluated at $t=\pm 1$.

The steepest-descent integral is optimal since it minimizes oscillations and ensures fast (exponential) decay. Nevertheless, the numerical integrations may be spoiled when $\theta_0$ is close to $\pm\pi/2$ because the integration path crosses over a logarithmic singularity, slowing down the convergence. Therefore the integrations can be regularized by updating $\theta_0$ as $\theta_0[1-\delta]$, where $\delta$ is small as shown in numerical results (Fig.~\ref{fig:result_1}).

\begin{figure*}[t]
    \centering\includegraphics[width=1\textwidth]{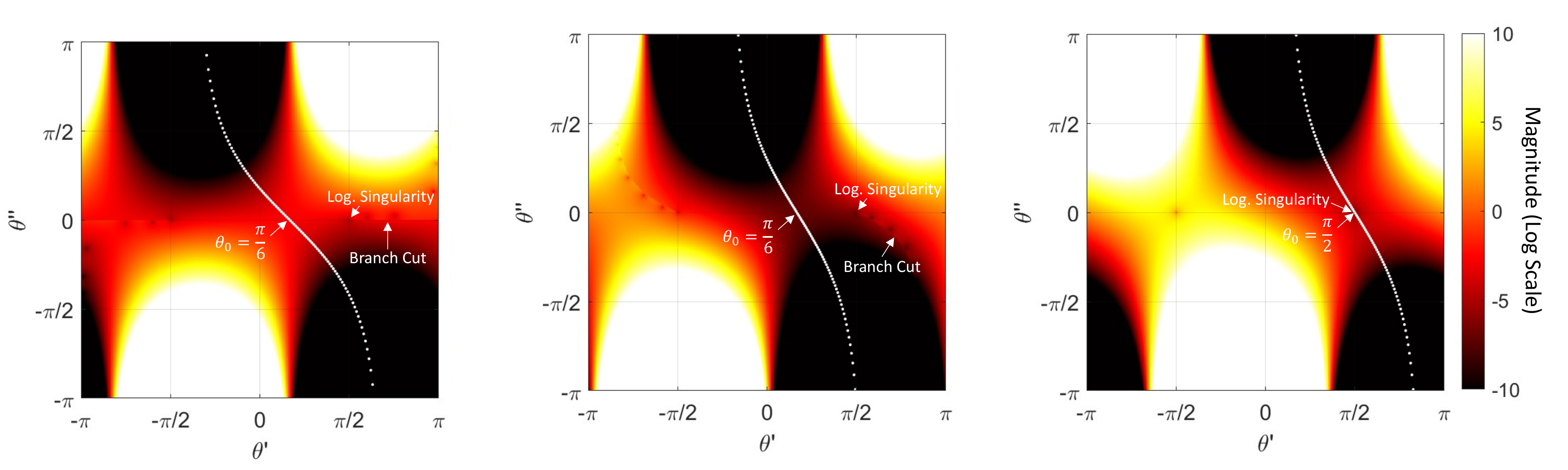}
    \footnotesize (a) \hspace{5.4cm} (b) \hspace{5.2cm}  (c)
  \caption{2.5-D Green's functions in the $\theta$ plane and the integration points on the steepest-descent contour with 100-point Gauss-Legendre quadrature. a) Lossless case with $r=\sqrt{2}\lambda$ and $\theta_0=\pi/6$, b) Lossy case with $k_0=2\pi[1+0.5i]$, c), Lossy case with $\theta_0=\pi/2$.}\label{fig:complex_plots}
\end{figure*}

\subsection{Steepest Descent Integration on the $s$ Plane}
Previous work~\cite{chew1990waves} parametrizes the steepest descent on the $s$ plane by making a transformation
\begin{equation}\label{eq:change_s}
    -s^2=f(\theta)-f(\theta_0).
\end{equation}
In this case, the real line in the $s$ plane maps to the steepest-descent contour in the $\theta$ plane because
\begin{equation}
    \Re\{f(\theta_0)\}\ge\Re\{f(\theta)\}
\end{equation}
on real line on the $s$ plane, and $\theta$ follows the constant-phase path as shown in (\ref{eq:constant_phase}). As a result, we can rewrite (\ref{eq:parts}) as
\begin{equation}
    g(r) = e^{ik_0r}\frac{i}{8\pi}\int_{\Gamma_\theta} d\theta\,\widehat{H}(\theta) e^{-s^2}.
\end{equation}
We can parametrize $s\in(-\infty,\infty)$ to the the $\theta$ plane as
\begin{equation}\label{eq:mapping_s}
    \theta = \theta_0+\cos^{-1}\left(1-\frac{s^2}{ik_0r}\right)\operatorname{sign}(s),
\end{equation}
where $\operatorname{sign}$ is the signum function that is used to choose the physical branch in the $s$ plane.

With the mapping in (\ref{eq:mapping_s}), we can write the correction term due to the change of variable, i.e., $d\theta=c(s)ds$, as
\begin{equation}
    c(s) = -\frac{2i}{\sqrt{s^2-2ik_0r}}.
\end{equation}
As a result, we can compute the integral on the real line as
\begin{equation}
    g(\boldsymbol{r},\boldsymbol{r}') = e^{f(\theta_0)}\frac{1}{2\pi}\int_{-\infty}^\infty ds\,e^{-s^2}c(s)\widehat{H}(\theta).
\end{equation}

The infinite integral can be evaluated numerically using Gauss-Hermite (G-H) quadrature~\cite{abramowitz1988handbook} as
\begin{equation}
    I=\int_{-\infty}^\infty e^{-x^2} h(x)\,dx \approx \sum_{i=1}^n w_i h(x_i)
\end{equation}
with $n$ points, where
\begin{equation}
    w_i=\frac{2^{n-1}n!\sqrt{\pi}}{n^2[H_{n-1}(x_i)]^2}
\end{equation}
is the quadrature weight and $H_{n}$ is the $n$\textsuperscript{th}-order Hermite polynomial.

\subsection{Steepest Descent Integration on the t Plane}

We can also map the steepest-descent contour on a finite interval $t\in(-\pi/2,\pi/2)$ by a transformation
\begin{equation}
    -\tan^2 t = f(\theta)-f(\theta_0).
\end{equation}
As a result, we can rewrite the mapping of of the SD countour as
\begin{equation}
    \theta = \theta_0+\operatorname{sgn}(t)\cos^{-1}\left(1-\frac{\tan^2 t}{ik_0r}\right)
\end{equation}
and convert (\ref{eq:parts}) as
\begin{equation}
    g(r)=e^{f(\theta_0)}\frac{i}{8\pi}\int_{-\pi/2}^{\pi/2} dt\,c(t)\widehat{H}(\theta)e^{-\tan^2 t},
\end{equation}
where
\begin{equation}
    c(t) = -\frac{2i\sec^2 t}{\sqrt{\tan^2t-2ik_0r}}.
\end{equation}
The integral above can be evaluated with G-L quadrature in the form of:
\begin{equation}
    \int_{-1}^1 h(x)dx\approx\sum_{i=1}^nw_ih(x_i)
\end{equation}
with a change of interval as similar to (\ref{eq:3greens_function_eff}).

\begin{figure*}[h!]
    \centering\includegraphics[width=0.9\textwidth]{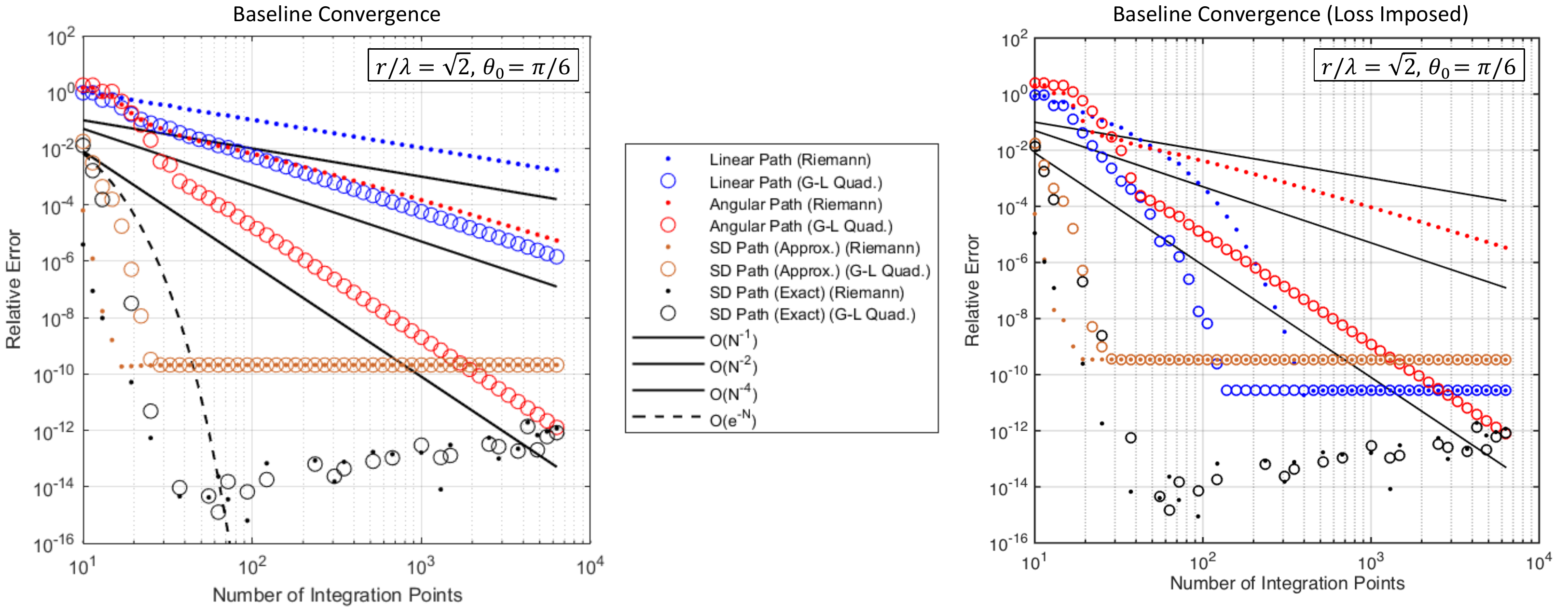}\\
  \vspace{0mm}{\hspace{0.4cm}\small(a)\hspace{9.5cm}(b)}
  \caption{Convergence with Riemann and Gaussian-Legendre integration rules when $r=\sqrt{2}\lambda$, $\theta_0=\pi/6$ for (a) lossless medium and (b) with imposed loss as regularizer.}\label{fig:result_1}
\end{figure*}

\section{Numerical Results}
In order to demonstrate and compare convergence performances of proposed integration paths in Sections~\ref{sec:quadratic_approx} and~\ref{sec:exact_SD}, several scenarios in lossless medium are considered. In all experiments, a point source is located at the origin and the radiation is observed from distance $r$ and observation angle $\theta_0$ as depicted in Fig.~\ref{fig:observation}. Here, $r$ is given in wavelength ($\lambda$) in free space. In order to improve convergence rates, several regularization techniques are proposed as summarized in Table \ref{tab:summary}. Effects of the proposed regularization techniques are discussed for each case in detail.

\begin{table}[h!]
\centering
\caption{Summary of Experimental Cases and Regularization Techniques}
\label{tab:summary}
\begin{tabular}{r|r|l}
\textbf{Figure} & \textbf{Observation}                     & \textbf{Regularization} \\ \hline
Fig.~\ref{fig:result_1}. a)     & $r=\sqrt{2}\lambda$, $\theta_0=\pi/6$    & -                      \\
Fig.~\ref{fig:result_1}. b)     & $r=\sqrt{2}\lambda$, $\theta_0=\pi/6$    & Loss Imposed            \\
Fig.~\ref{fig:result_2}. a)     & $r=\sqrt{2}\lambda$, $\theta_0=\pi/2$    & -                      \\
Fig.~\ref{fig:result_2}. b)     & $r=\sqrt{2}\lambda$, $\theta_0=\pi/2$    & Path Shifted           \\
Fig.~\ref{fig:result_2}. c)     & $r=0.1\sqrt{2}\lambda$, $\theta_0=\pi/6$ & -                      \\
Fig.~\ref{fig:result_3}. a)     & $r=0.1\sqrt{2}\lambda$, $\theta_0=\pi/6$ & Limits Increased        \\
Fig.~\ref{fig:result_3}. b)     & $r=10\sqrt{2}\lambda$, $\theta_0=\pi/6$  & -              
\end{tabular}
\end{table}

\begin{figure*}[h!]
    \centering\includegraphics[width=1\textwidth]{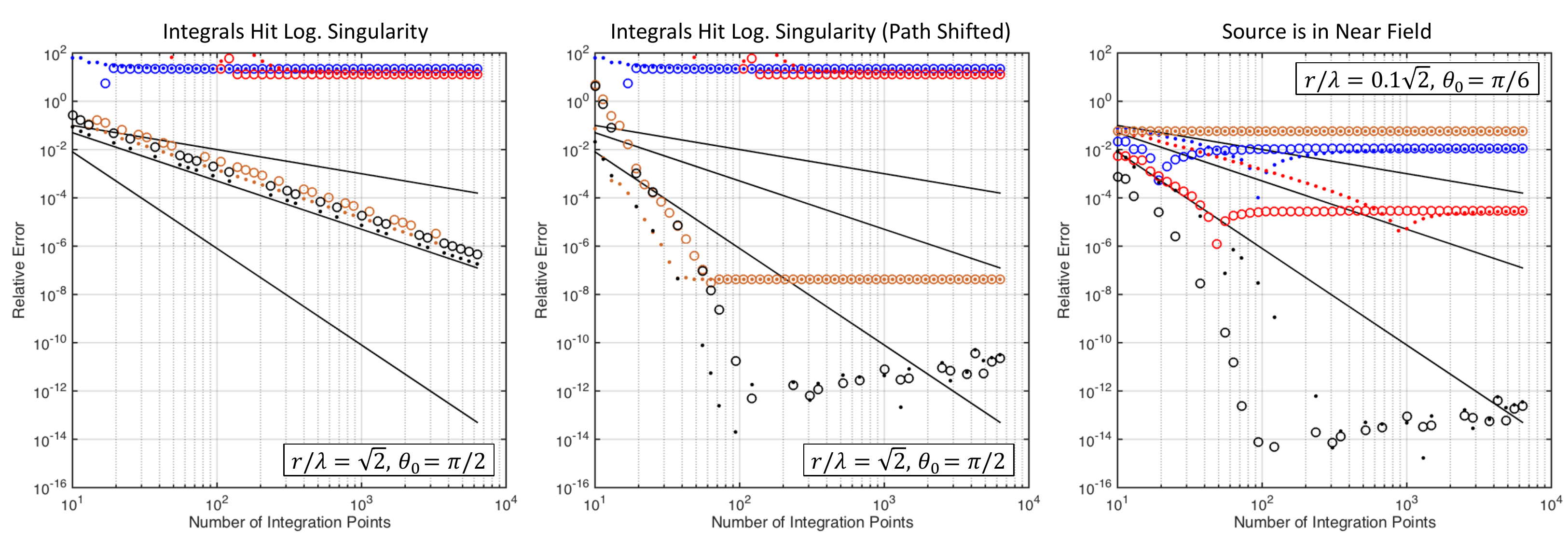}\\
  \vspace{-3mm}{\hspace{0.5cm}\small(a)\hspace{5.8cm}(b)\hspace{5.7cm}(c)}
  \caption{Convergence with Riemann and Gaussian-Legendre integration rules when $r=\sqrt{2}\lambda$, $\theta_0=\pi/2$ with (a) original convergence when integration contours hit logarithmic singularity and (b) regularized convergence by shifting the integration contours. (c) The observation is near to the source as $r=0.1\sqrt{2}\lambda$, $\theta_0=\pi/6$ where truncated integrals do not capture evanescent the spectrum and thus do not converge (except the exact steepest-descent integration).}\label{fig:result_2}
\end{figure*}

\begin{figure*}[h!]
    \centering\includegraphics[width=1\textwidth]{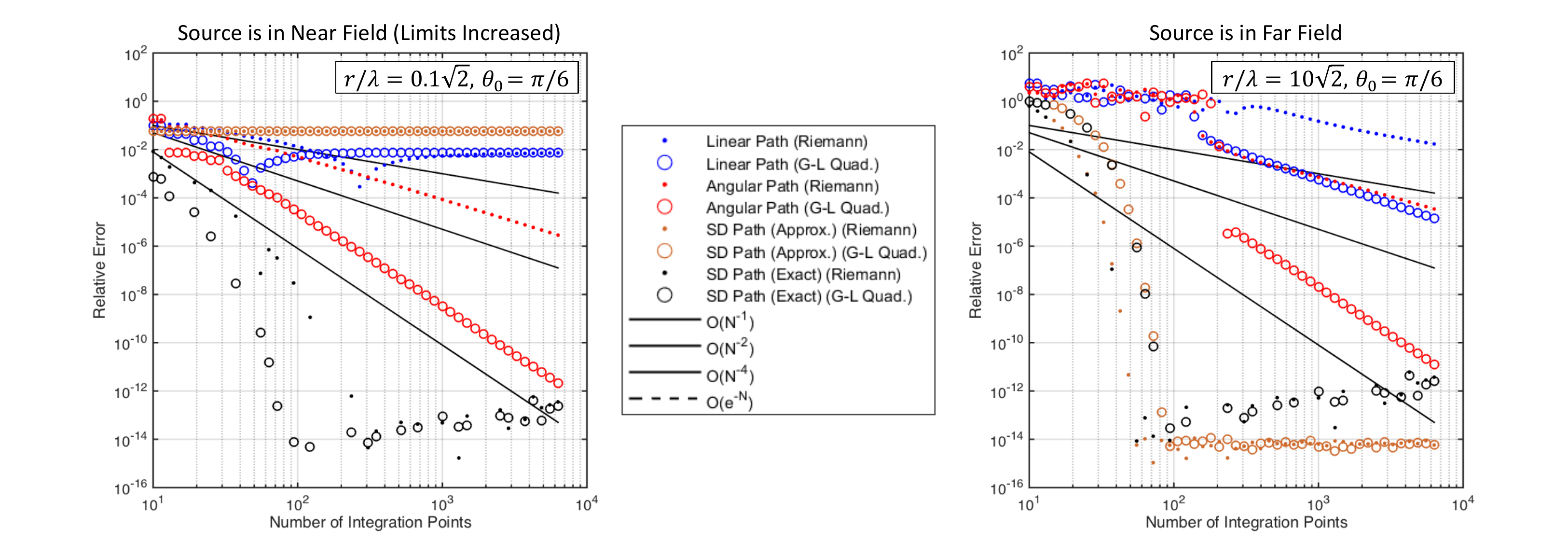}\\
  \vspace{-3mm}{\hspace{0.4cm}\small(a)\hspace{10cm}(b)}
  \caption{Convergence with Riemann and Gaussian-Legendre integration rules when $r=0.1\sqrt{2}\lambda$, $\theta_0=\pi/6$ with a) increased integration limits. (b) When the observation is far away to the source as $r=10\sqrt{2}\lambda$, $\theta_0=\pi/6$.}\label{fig:result_3}
  \vspace{-2mm}
\end{figure*}

Figures~\ref{fig:result_1}--\ref{fig:result_3} show relative errors of numerical integrations with respect to the number of integration points, where point markings show the Riemann integration error whereas circle markings show Gaussian-Legendere (denoted as G-L) integration results. The relative error is defined as
\begin{equation}
    E = r|e^{ik_0r}-4\pi I|,
\end{equation}
where $I$ is the numerical integration result. The figures involve solid lines showing $\mathcal{O}(N^{-1})$, $\mathcal{O}(N^{-2})$, and $\mathcal{O}(N^{-4})$ curves as reference to linear, quadratic, and fourth-order convergence rates, respectively, where $N$ is the number of integration points. Results are colored as the following.

The blue data points use linear, i.e., $k_z$, integration in (\ref{eq:hankel}) truncated as $k_z\in(-2k_0, 2k_0)$. This partially captures the evanescent-wave spectrum and fully captures the propagating-wave spectrum, as depicted in Fig.~\ref{fig:kz_plane}. The red data points use angular integration in (\ref{eq:angular}), whose contour is depicted in Fig.~\ref{fig:kt_plane}. The imaginary angle $\theta''$ is truncated as $\theta''\in(-1.5\pi, 1.5\pi)$. As a result of this truncation, this integration partially captures the evanescent-wave spectrum. The brown data points use the quadratic approximation to the steepest-descent path, as proposed in (\ref{eq:approximate_descent}). Finally, the black data points use the steepest-descent path as proposed in (\ref{eq:steepest_descent}), which captures the full-wave spectrum. The steepest-descent contour is depicted in Fig.~\ref{fig:kt_plane_efficient}.


Fig.~\ref{fig:result_1} (a) shows convergence of integration when $r=\sqrt{2}\lambda$, $\theta_0=\pi/6$. Evidently, the linear integration has linear and quadratic convergence with Riemann and Gauss-Legendre integrations, respectively. Similarly, the angular integration has quadratic and fourth-order convergence with the respective integration rules. The poor performance of the linear integration compared to the angular spectrum path is a result of double and single branch points with which they both coincide, respectively. In other words, the linear contour coincides with algebraic and logarithmic branch points, whereas the angular contour path coincides only with logarithmic branch points (and the angular spectrum formulation removes the algebraic branch points). In contrast, approximate and exact steepest-descent integrations have exponential convergence because they avoid the logarithmic branch points. Nevertheless, the convergence of the approximate integration is limited. Fig.~\ref{fig:result_1} (b) shows the case when a slightly small loss is imposed ($k_0''=0.05$, where $k_0=k_0'+ik_0''$) as a regularizer to accelerate convergence. In this case, the linear path enjoys exponential convergence because the branch points move away from the integration path as seen in Fig.~\ref{fig:kz_plane}. However, other integrations are not effected because the loss does not move the branch points in those cases. Also, we can observe that the linear integration converges to a less accurate result because of regularization. It is worthy to comment that no path (except the steepest-descent path) can capture full spectrum because of truncations in these experiments, but their range are enough to observe the convergence trends up to double precision.

Fig.~\ref{fig:result_2} (a) shows when the observation angle $\theta_0=\pi/2$. This case is problematic because the integration contours move closer to the logarithmic singularity the complex plane, and therefore the their convergence are hampered. In order to improve the convergence rates, we propose shifting the integration paths such that they passes through $\theta_0[1-\delta]$, where $\delta$ is small. Fig.~\ref{fig:result_2}(b) shows exponential convergence for the steepest-descent integrations when $\delta=0.1$. As a result of the regularization, the accuracy of the final results are slightly reduced, however, it is rewarded with the improved convergence.

Fig.~\ref{fig:result_2}(c) shows the convergence when the observation point is close to the source, where $r=0.1\sqrt{2}$, and $\theta_0=\pi/6$. In this case, only the exact steepest-descent path can capture the full spectrum because of the truncation in other cases. As a result, their convergence are saturated. For improving the convergence, we can increase the integral limits in order to capture more of the evanescent spectrum: Fig.~\ref{fig:result_3}(a) shows the convergence results when the integral limits of the truncated integrals are doubled to capture the more of the evanescent spectrum. The result shows that the angular spectrum integration has significant accuracy improvement whereas the linear integration has a slight improvement. This result can be explained by the exponential decay rates of the evanescent spectrum with the angular spectrum formulation. As a result, with the increased integration limits, the angular spectrum integration regain its convergence rates with Riemann and Gaussian-Legendre quadrature rules.

Lastly, we place the observation far away from the source, where the approximate steepest-descent converges down to double precision (accurate up to 14 significant decimal digits) within a few integration points. The exact path does not attain double precision because of the numerical approximations in computing trancendental functions, such as the Hankel function in Eq.~\ref{eq:steepest_descent}.

\section{Conclusion}

In this paper, we analyze the expansion of the 3-D Green's function with 2-D cylindrical modes in the context of solving 2.5-D inverse problems. We discuss the logarithmic singularity and its suppression with the angular-spectrum formulation. Furthermore, propose a novel parametrization for the exact steepest-descent for numerical implementation that is suitable to solve 2.5-D inverse scattering problems with lower computational complexity. Alternatively, we propose an approximate (yet simple) steepest-descent integration technique when the observation is far from the source (in terms of wavelengths). Numerical results show that the proposed integration techniques have exponential convergence, meaning that requires less number of 2-D solutions to solve a 2.5-D inverse problem.

\section*{Acknowledgments}
This work was supported in part by grants from the NIH (R21EB023403 and R01CA251939). The authors greatly acknowledge John Clark's help on drawing three-dimensional figures.

\bibliographystyle{ieeetr}
\bibliography{ref}

\newpage

\section{Appendix: Sommerfeld Identities for 2.5-D}\label{sec:appendix}

\begin{figure}[t]
  \centering
  \hspace{-4mm}
 \includegraphics[width=0.5\textwidth]{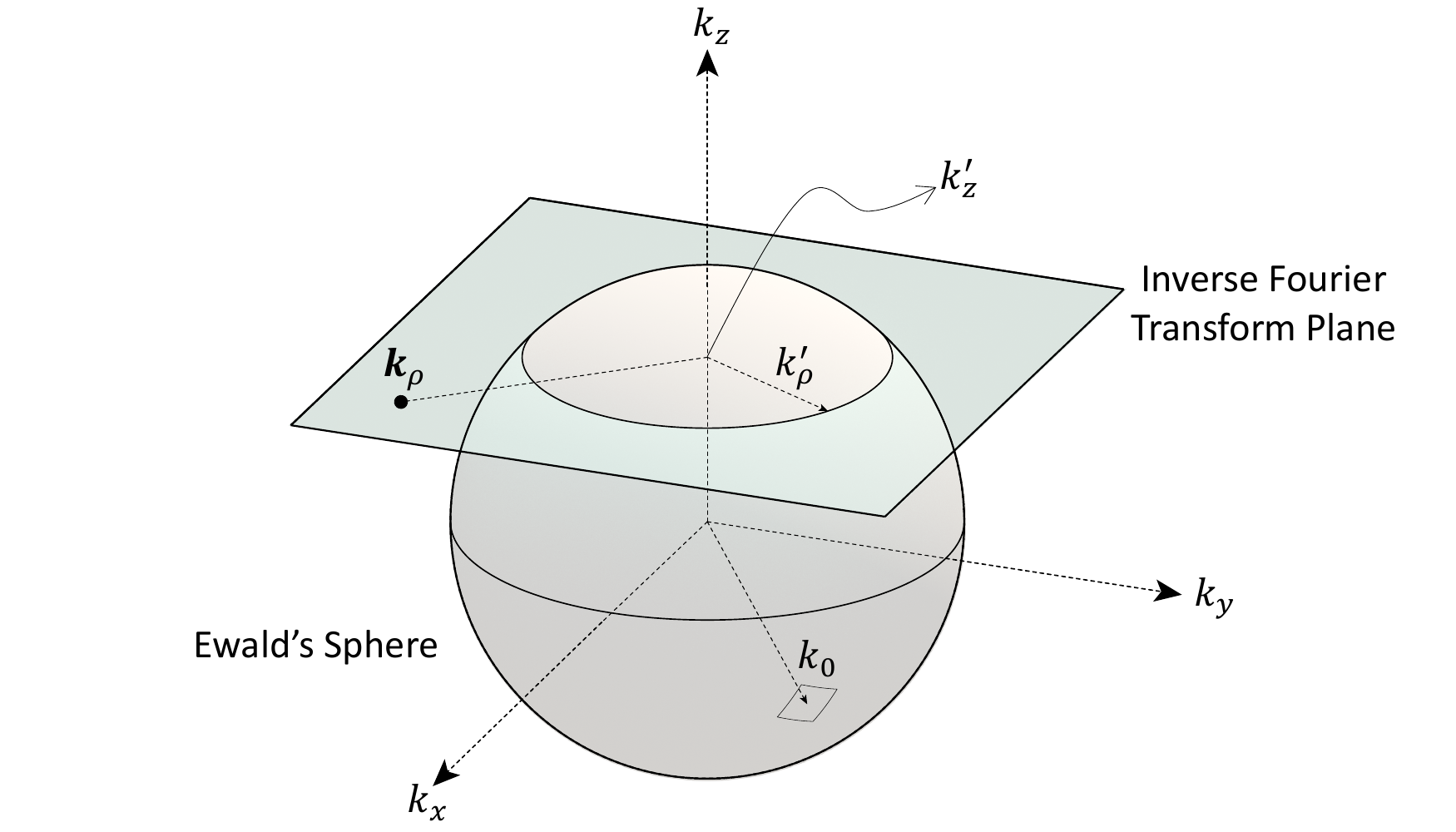}
\vspace{-3mm}
  \caption{Ewald's sphere and 2-D inverse Fourier transform plane at $k_z'$.}\label{fig:ewald_sphere}
\end{figure}

In order to find a closed-form expression of the 2.5-D Green's function in (\ref{eq:2p5_green}), we can take a spatial Fourier transform of the wave equation in (\ref{eq:wave_freq}) as
\begin{equation}
    \int d\boldsymbol{r}\, e^{-i\boldsymbol{k}\cdot\boldsymbol{r}}[\nabla^2g_0(\boldsymbol{r},\boldsymbol{r}')+k_0^2g_0(\boldsymbol{r},\boldsymbol{r}')=-\delta(\boldsymbol{r}-\boldsymbol{r}')]
\end{equation}
yielding a closed-form expression in Fourier domain, i.e.,
\begin{equation}
    g_0(\boldsymbol{k},\boldsymbol{r}') = \frac{e^{-i\boldsymbol{k}\cdot\boldsymbol{r}'}}{k^2-k_0^2}.
\end{equation}
The above follows from the fact that $\nabla^2 e^{-i\boldsymbol{k}\cdot\boldsymbol{r}}=-k^2e^{-i\boldsymbol{k}\cdot\boldsymbol{r}}$. Then we let $k_0^2=k_\rho'^2+k_z'^2$ and perform a partial inverse Fourier transform on the $k_\rho$ plane  as
\begin{align}
    \begin{split}
        g_0(\boldsymbol{\rho},\boldsymbol{r}') &=\frac{1}{(2\pi)^2}\int d\boldsymbol{k}_\rho\, e^{i\boldsymbol{k}_\rho\cdot\boldsymbol{\rho}}g_0(\boldsymbol{k},\boldsymbol{r}')\\
        &=\frac{e^{-ik_zz'}}{(2\pi)^2}\int d\boldsymbol{k}_\rho \frac{e^{i\boldsymbol{k}_\rho\cdot(\boldsymbol{\rho}-\boldsymbol{\rho}')}}{k_\rho^2-k_\rho'^2}\label{eq:partial_fourier}.
    \end{split}
\end{align}
Note that $k_z'$ is fixed, and $k_\rho'$ is constrained on the Ewald's sphere depicted in Fig.~\ref{fig:ewald_sphere}. This is defined as the 2.5-D free-space Green's function and is denoted as $g_0(\boldsymbol{\rho},\boldsymbol{r}')$ in order to differentiate it from the 3-D Green's function $g_0(\boldsymbol{r},\boldsymbol{r}')$ and the 2-D Green's function $g_0(\boldsymbol{\rho},\boldsymbol{\rho}')$.

To find the closed-form of the 2.5-D Green's function, we need to perform the inverse Fourier integral in (\ref{eq:partial_fourier}). This subsection assumes that the point source is located at the origin of the cylindrical coordinate system, i.e., $\boldsymbol{r}'=z'\hat{\boldsymbol{z}}$, i.e., $\boldsymbol{\rho}'=0$ for brevity. Using the cylindrical symmetry, we can set $\boldsymbol{\rho}=\rho\hat{\boldsymbol{x}}$ and also rewrire $\boldsymbol{\rho}=\hat{\boldsymbol{x}}\rho\cos{\alpha}+\hat{\boldsymbol{y}}\rho\sin{\alpha}$ and $d\boldsymbol{k}_\rho=dk_xdk_y=k_\rho dk_\rho \alpha$ in polar form. Then we can write the integral (\ref{eq:partial_fourier}) as
\begin{equation}\label{eq:bessel}
    g_0(\rho,z') = \frac{e^{-ik_zz'}}{2\pi}\int_0^\infty dk_\rho \frac{k_\rho}{{k_\rho^2-k_\rho'^2}}\underbrace{\frac{1}{2\pi}\int_{-\pi}^{\pi}  e^{ik_\rho\rho\cos\alpha} d\alpha}_{J_0(k_\rho\rho)}.
\end{equation}
The inner part (noted with curly braces) is the definition~\cite{abramowitz1988handbook} of the zeroth-order Bessel function, denoted as $J_0(z)$. Consequently, an exhaustive list of identities can be written as
\begin{align}
    g_0(\rho,z') 
    &= \frac{e^{-ik_z'z'}}{2\pi}\int_0^\infty  dk_\rho\frac{k_\rho}{k_\rho^2-k_\rho'^2}J_0(k_\rho\rho)\\
    &= \frac{e^{-ik_z'z'}}{4\pi}\int_{-\infty}^\infty  dk_\rho\frac{k_\rho }{k_\rho^2-k_\rho'^2}J_0(k_\rho\rho)\\
    \begin{split}
        =\frac{e^{-ik_z'z'}}{4\pi}\left[\int_0^\infty dk_\rho\frac{k_\rho }{k_\rho^2-k_\rho'^2}H_0^{(1)}(k_\rho\rho)\right.\\
        +\left.\int_0^\infty dk_\rho\frac{k_\rho }{k_\rho^2-k_\rho'^2}H_0^{(2)}(k_\rho\rho)\right]
    \end{split}&\\
    \begin{split}
        =\frac{e^{-ik_z'z'}}{8\pi}\left[\int_{-\infty}^\infty dk_\rho\frac{k_\rho }{k_\rho^2-k_\rho'^2}H_0^{(1)}(k_\rho\rho)\right.\\
        +\left.\int_{-\infty}^\infty dk_\rho\frac{k_\rho }{k_\rho^2-k_\rho'^2}H_0^{(2)}(k_\rho\rho)\right]
    \end{split}&\\
    &= \frac{e^{-ik_z'z'}}{4\pi}\int_{-\infty}^\infty  dk_\rho\frac{k_\rho }{k_\rho^2-k_\rho'^2}H_0^{(1)}(k_\rho\rho)\label{eq:sommerfeld_last}\\
    &= \frac{e^{-ik_z'z'}}{4\pi}\int_{-\infty}^\infty  dk_\rho\frac{k_\rho }{k_\rho^2-k_\rho'^2}H_0^{(2)}(k_\rho\rho)
\end{align}
which can be regarded as Sommerfeld identities for 2.5-D\footnote{Further reading on \S\S2.2.2 of \cite{red} is encouraged.}.

\begin{figure}[t]
    \centering\includegraphics[width=0.43\textwidth]{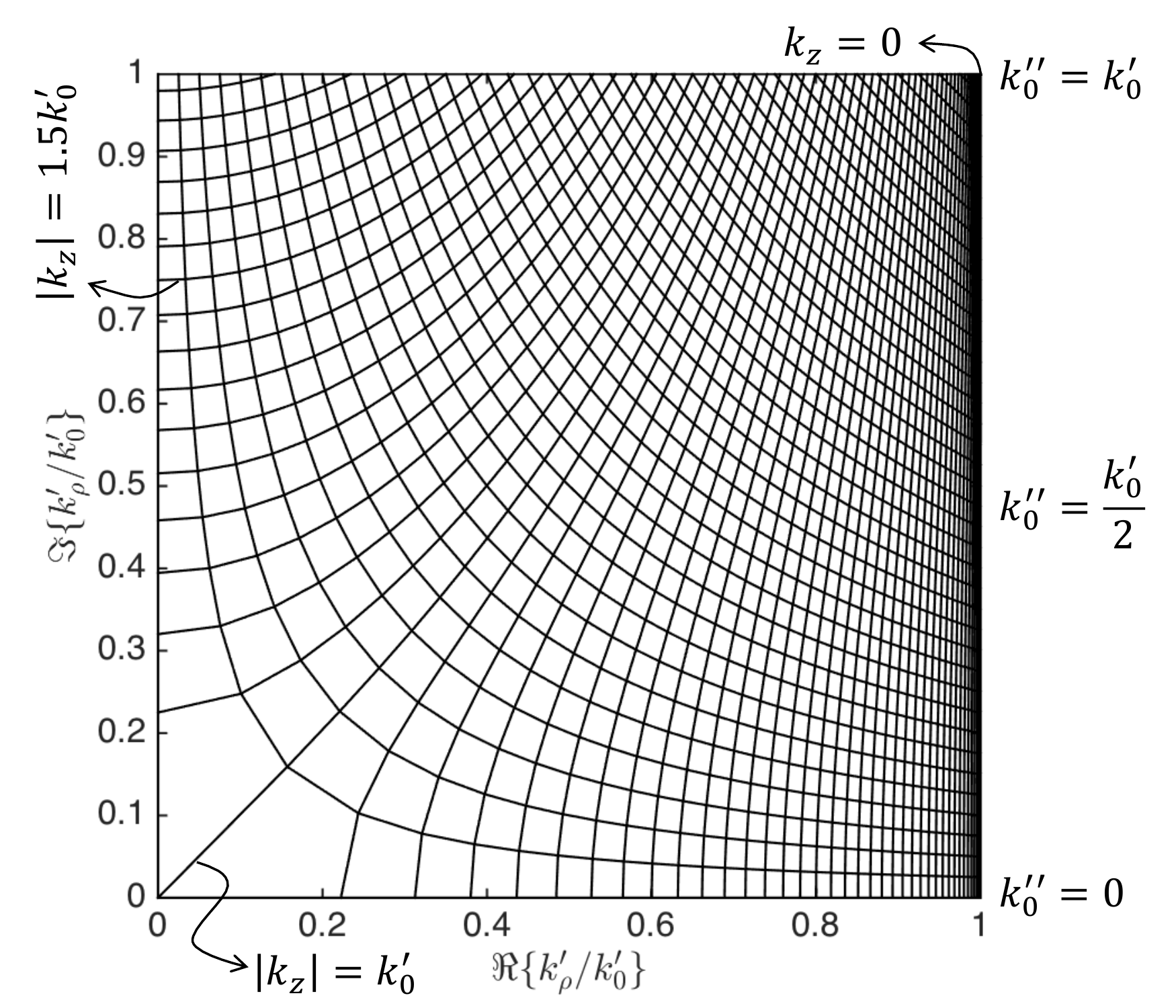}
  \vspace{-3mm}
  \caption{Conformal mapping of $k_\rho$ for various real $k_z$ and complex $k_0$ values with $k_0'/40$ step size. For increasing $k_z$ and fixed $k_0$, $k_\rho$ travels left-to-right and right-to-left (on a constant $k_0$-line) when $k_z<0$ and $k_z>0$, respectively. For fixed $k_z$ and increasing $k_0''$, $k_\rho$ travels upwards (on a constant $k_z$-line).}\label{fig:grid}
\end{figure}

These integrals can be seen as the cylindrical Fourier expansion of the Green's function with different kind of modes, which can be physically interpreted as standing- and propagating-wave modes. The first expansion uses standing-wave modes described by Bessel functions. We can extend the first integral to negative frequencies by an unfolding using a symmetry in (\ref{eq:bessel}), i.e., we can limit the inner integral between $(0,\pi)$ and let the outer limits be $(-\infty,\infty)$, which leads to the second identity. The latter two integrals write the standing waves as a sum of incoming and outgoing propagating waves, i.e., $2J_0(z) = H_0^{(1)}(z)+H_0^{(2)}(z)$, described by zeroth-order Hankel functions of the first and second kind, respectively. And the last two follows from the connection formula $H_0^{(1)}(z)=e^{i\pi}H_0^{(2)}(e^{i\pi}z)$ and a change of variable.

\begin{figure}[h!]
    \vspace{-3mm}
  \centering\includegraphics[width=0.48\textwidth]{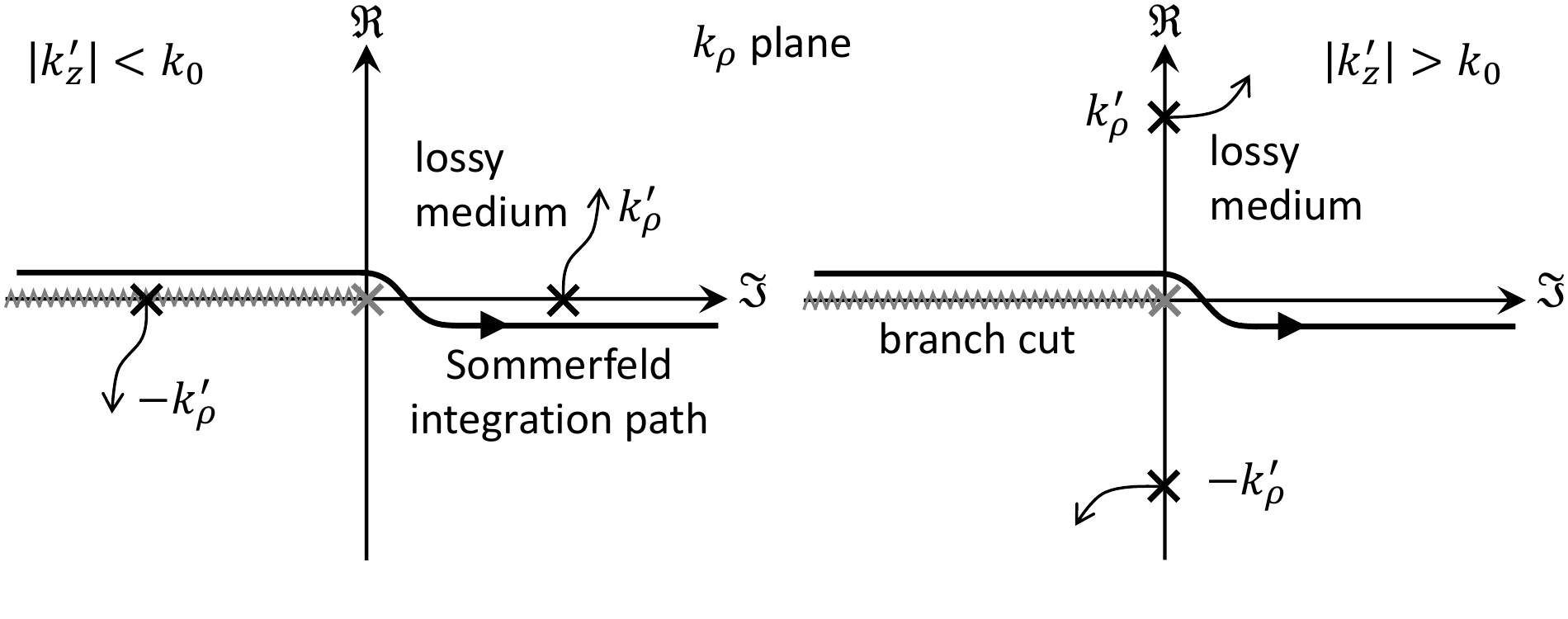}
  \vspace{-3mm}
  \caption{The $k_\rho$-plane for the integrand of (\ref{eq:sommerfeld_last}). The Sommerfeld integration path picks up the physical pole on the first quadrant of the $k_\rho$-plane to satisfies the Sommerfeld radiation condition.}\label{fig:sommerfeld_integration}
\end{figure}

We can perform (\ref{eq:sommerfeld_last}) analytically with the residue theorem by closing the integration path in the upper-half of the complex plane, where the integrand vanishes as $|k_\rho|\rightarrow\infty$. The integrand involves a logarithmic branch cut as a result of the Hankel function, and simple poles at $\pm k_\rho'$. For the lossless case, i.e., when $k_0$ is purely real, the poles are depicted in Fig.~\ref{fig:sommerfeld_integration}. In this case, the poles $\pm k_\rho'$ are purely real when $|k_z'|<k_0$ and purely imaginary when $|k_z'|>k_0$, respectively. The former case is more problematic because the result may violate the Sommerfeld radiation condition, depending on which poles we pick.

To pick the physical pole, we impose a loss to the medium by letting $k_0 = k_0'+ik_0''$, where $k_0''>0$, which shifts the poles slightly out of the axes. The direction of the shifts are depicted with arrows in Fig.~\ref{fig:sommerfeld_integration}.  The pole on the negative real axis shifts below of the original integration contour (to the third quadrant), and the pole on the positive axis shifts above the integration contour (to the first quadrant). A set of loci of $k_\rho'$ for different values of $k_0$ and $k_z'$ is shown in Fig.~\ref{fig:grid} (with unprimed coordinates). Evidently, the pole on the first quadrant is the physical pole. By analytic continuation, we can deform the integration path as shown in Fig.~\ref{fig:sommerfeld_integration} to avoiding the logarithmic branch and pick the physical pole. This path is known as the Sommerfeld integration path. The lossless case is not physical anyways, since the energy from a point source should vanish at infinity, which is known as the Sommerfeld radiation condition~\cite{schot1992eighty}. To satisfy, we consider the lossless medium as the limiting case of lossy medium, i.e., $k_z''\rightarrow0$ to satisfy this condition. Invoking the residue theorem, we can write the Sommerfeld identity
\begin{equation}
    \int_{-\infty}^\infty  dk_\rho\frac{k_\rho }{k_\rho^2-k_\rho'^2}H_0^{(1)}(k_\rho\rho)= i\pi H_0^{(1)}(k_\rho'\rho).
\end{equation}
As a result, switching the unprimed coordinates to denote the spectral points on the Ewald's spehere,
\begin{equation}
    g_0(\rho,z') = \frac{i}{4}H_0^{(1)}(k_\rho\rho)e^{-ik_zz'}.\label{eq:sommerfeld_identity}
\end{equation}
Likewise, the other Sommerfeld identities can be written similar to the above. Substituting this identity for a point source located at an arbitrary location $\boldsymbol{r}'$, we arrive (\ref{eq:2p5_green}).

\section{Appendix: Derivation of 3-D Green's Function}
For analytical integration of (\ref{eq:hankel}), we can can truncate its interval and rotate the cylindrical coordinate system by letting $z\rightarrow r$ (equivalently by letting $\rho\rightarrow0$ as
\begin{multline}
    \lim_{z\rightarrow r} \frac{i}{8\pi}\int_{-A}^{A} dk_z H_0^{(1)}(\sqrt{k_0{}^2-k_z^2}\,\rho)e^{ik_zz}\\
    = 
    \lim_{z\rightarrow r} \frac{i}{8\pi}\int_{-A}^{A} dk_z \frac{2i}{\pi}\left(\log\frac{\sqrt{k_0{}^2-k_z^2}\,\rho}{2}+\gamma\right)e^{ik_zz},
\end{multline}
by using the asymptotic form of Hankel function for small argument, where $\gamma$ is the Euler-Mascheroni number~\cite{abramowitz1988handbook}. Integrating the above by parts,
\begin{multline}\label{eq:integrate_parts}
    \lim_{z\rightarrow r} \frac{-1}{4\pi^2} \left[ \left(\log\frac{\sqrt{k_0{}^2-k_z^2}\,\rho}{2}+\gamma\right)  \frac{e^{ik_zz}}{ir} \right]_{-A}^{+A}\\
    + \frac{1}{4\pi^2}\int_{-A}^{+A} dk_z \frac{k_z}{k_z^2-k_0^2}\frac{e^{ik_zr}}{ir}.
\end{multline}
Letting $A\rightarrow\infty$, one can show that the first term vanishes and the second term picks up the pole at $k_0$, yielding $g_0(r) = e^{ik_0r}/{4\pi r}$.

\begin{figure}[!ht]
  \centering\includegraphics[width=0.45\textwidth]{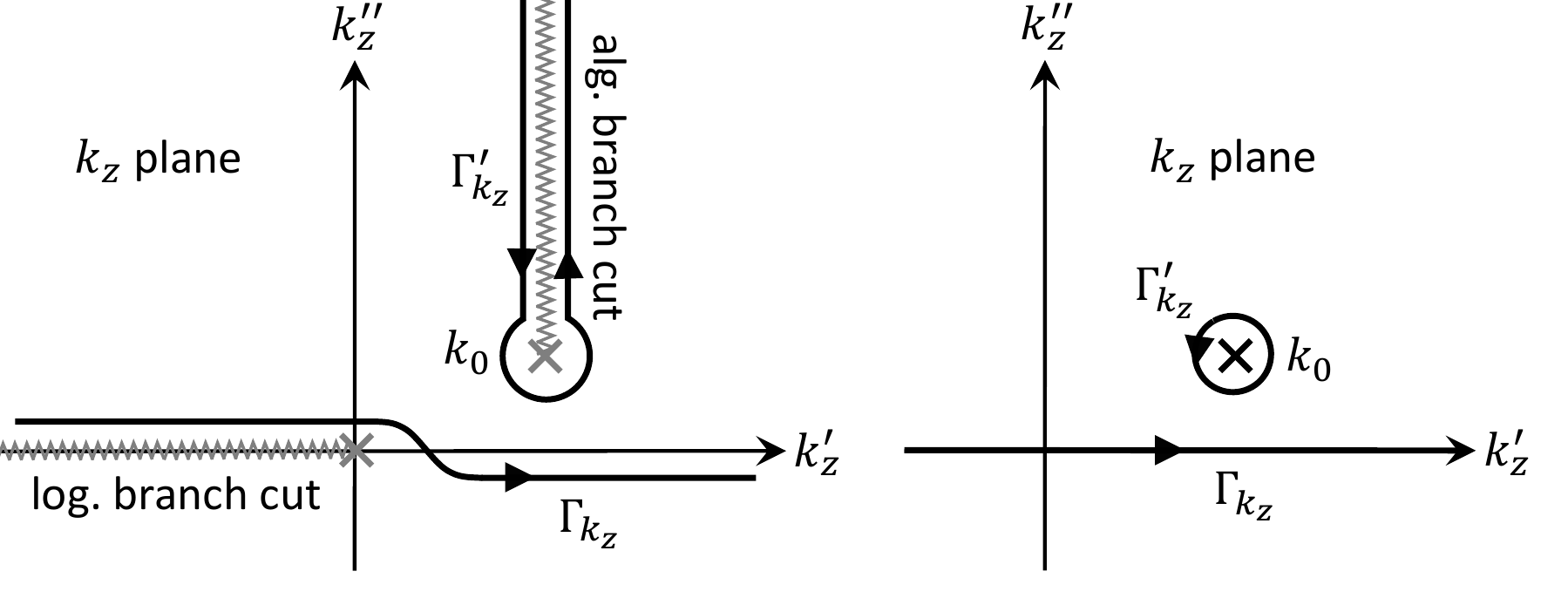}\\
  {\small(a)\hspace{3.5cm}(b)}
  \caption{(a) The modified branch cuts and integration paths for the first term (a) and the second term (b) of (\ref{eq:integrate_parts}).}\label{fig:kz_analytical}
\end{figure}

Specifically, invoking the Jordan's lemma and analytical continuation, the integral path algebraic branch cut connected to $k_0$ can be chosen to go $k_0'+i\infty$. Then the first term vanishes since $\log(\sqrt{k_0{}^2-k_z^2}\rho)e^{ik_zz}\rightarrow0$, when $k_z\rightarrow i\infty$
Similarly, it is easy to see that the first term is free of branch points and picks up the pole at $k_0$. The composite $k_z$-plane involving modified branch cuts and integration paths for the first and terms of (\ref{eq:integrate_parts})  is shown in Fig.~\ref{fig:kz_analytical}(a) and Fig.~\ref{fig:kz_analytical}(b), respectively.

This result confirms (\ref{eq:fourier_synthesis}) as
\begin{equation}
    \frac{i}{8\pi}\int dk_z H_0^{(1)}(k_\rho'|\boldsymbol{\rho}-\boldsymbol{\rho}'|)e^{ik_z(z-z')}=\frac{e^{ik_0|\boldsymbol{r}-\boldsymbol{r}'|}}{4\pi|\boldsymbol{r}-\boldsymbol{r}'|}
\end{equation}
with the point source placed at arbitrary location. The right-hand-side of this identity can also be obtained by another method, such as singularity matching. In any case, since both sides satisfies the Sommerfeld radiation condition~\cite{schot1992eighty}, by the unqiueness theorem, they have to be equal.

\end{document}